\pdfoutput=1
\documentclass[twoside]{article}
\usepackage[letterpaper, left=3cm, right=3cm, top=4cm, bottom=2.5cm]{geometry}
\usepackage{graphicx}
\usepackage{amsmath,amssymb,amsthm}
\usepackage{authblk}
\usepackage[utf8]{inputenc}
\usepackage{microtype}
\usepackage{mathrsfs} 
\usepackage{enumitem}
\usepackage{calc}
\usepackage{esint}
\usepackage{fancyhdr}
\usepackage{xcolor}
\usepackage[labelfont = bf]{caption}
\usepackage{doi}
\usepackage{subfigure}
\usepackage[numbers]{natbib}
\usepackage{float}
\usepackage{stmaryrd}
\usepackage{mathtools}
\usepackage{booktabs}
\usepackage{colortbl}
\usepackage{tabulary}
\usepackage{diagbox}
\usepackage{caption}
\usepackage{scalerel} 

\DeclareMathAlphabet{\mathbbold}{U}{bbold}{m}{n}

\definecolor{mordantred19}{rgb}{0.68, 0.05, 0.0}
\definecolor{pakistangreen}{rgb}{0.0, 0.4, 0.0}

\usepackage{physics} 
\usepackage{tikz}
\usepackage{mathdots}
\usepackage{yhmath}
\usepackage{cancel}
\usepackage{color}
\usepackage{siunitx}
\usepackage{array}
\usepackage{multirow}
\usepackage{amssymb}
\usepackage{gensymb}
\usepackage{tabularx}
\usepackage{extarrows}
\usepackage{booktabs}
\usetikzlibrary{fadings}
\usetikzlibrary{patterns}
\usetikzlibrary{shadows.blur}
\usetikzlibrary{shapes}

\newtheorem{theorem}{Theorem}[section]
\numberwithin{equation}{section}
\newtheorem{proposition}[theorem]{Proposition}
\newtheorem{definition}[theorem]{Definition}

\newtheorem{remark}[theorem]{Remark}
\newtheorem{lemma}[theorem]{Lemma}

\newtheorem{assumption}[theorem]{Assumption}

\usepackage{titlesec}
\titleformat{\section}{\normalfont\scshape\centering}{\thesection.}{0.5em}{}
\titleformat*{\subsection}{\itshape}
\titleformat*{\subsubsection}{\itshape}

\providecommand{\keywords}[1]
{
	{\small\emph{Keywords:} #1}
}

\providecommand{\MSC}[1]
{
	{\small\emph{AMS MSC (2020):  } #1}
}

\definecolor{denim}{rgb}{0.08, 0.38, 0.74}
\definecolor{byzantium}{rgb}{0.44, 0.16, 0.39} 
\definecolor{shamrockgreen}{rgb}{0.0, 0.62, 0.38} 

\usepackage{hyperref}
\hypersetup{
	colorlinks=true,
	linkcolor=denim,
	citecolor = shamrockgreen,
	filecolor=magenta, 
	urlcolor=byzantium,
}

\pdfoutput=1
\begin{document}
	\setlength{\abovedisplayskip}{5pt}
	\setlength{\belowdisplayskip}{5pt}
	\setlength{\abovedisplayshortskip}{5pt}
	\setlength{\belowdisplayshortskip}{5pt}

	\title{\vspace{-15mm}Finite element discretization of the steady, generalized Navier--Stokes equations for small shear stress exponents} 
	\author[1]{Alex Kaltenbach\thanks{Email: \url{kaltenbach@math.tu-berlin.de}}}
	\author[2]{Julius Je\ss{}berger\thanks{Email: \url{julius.jessberger@mathematik.uni-freiburg.de}}}
	\date{\today\vspace{-2.5mm}}
	\affil[1]{\small{Institute of Mathematics, Technical University of Berlin, Stra\ss e des 17.\ Juni 135, 10623 Berlin, GERMANY}}
    \affil[2]{\small{Department of Applied Mathematics, University of Freiburg, Ernst--Zermelo-Str. 1, 79104 Freiburg, GERMANY}\vspace{-2.5mm}}
	\maketitle

	\pagestyle{fancy}
	\fancyhf{}
	\fancyheadoffset{0cm}
	\addtolength{\headheight}{-0.25cm}
	\renewcommand{\headrulewidth}{0pt} 
	\renewcommand{\footrulewidth}{0pt}
	\fancyhead[CO]{\textsc{FE discretization of steady, generalized Navier-Stokes equations}}
	\fancyhead[CE]{\textsc{A. Kaltenbach and J. Je\ss berger}}
	\fancyhead[R]{\thepage}
	\fancyfoot[R]{}
	
	\begin{abstract}
		A finite element (FE) discretization for the  steady, incompressible, fully inhomogeneous, generalized Navier--Stokes equations is proposed. By the method of divergence  reconstruction operators, the formulation is valid for all shear stress exponents $p > \tfrac{2d}{d+2}$. The Dirichlet boundary condition is imposed strongly, using any discretization of the boundary data which converges at a sufficient rate.\linebreak
	\textit{A priori} error estimates for the velocity vector field and  kinematic pressure are derived and numerical experiments are conducted. These confirm the quasi-optimality of the 	\textit{a priori} error~estimate~for~the velocity vector field. 
	The 	\textit{a priori} error estimates  for the kinematic~\mbox{pressure}~are~\mbox{quasi-optimal}~if~${p \hspace{-0.1em}\leq\hspace{-0.1em} 2}$.
	\end{abstract}
	
	
	\MSC{35J60; 35Q35; 65N12; 65N15; 65N30; 76A05.}
	
	\section{Introduction}\thispagestyle{empty}\enlargethispage{15mm}

   \hspace*{5mm}The steady motion of a homogeneous, incompressible generalized Newtonian fluid can be modeled by the \emph{generalized Navier--Stokes equations}, \textit{i.e.},
\begin{align} \label{eq:fem:main_problem1}
	- \operatorname{div} \mathbf{S}(\mathbf{Dv}) + \operatorname{div}(\mathbf{v}\otimes\mathbf{v}) + \nabla q = \mathbf{f} \quad \textrm{ in } \Omega\,.
\end{align}
Here, $\Omega\subseteq \mathbb{R}^d$, $d\in \{2,3\}$, denotes either a bounded polygonal (if $d=2$) or~\mbox{polyhedral}~(if~${d=3}$) Lipschitz domain and $\mathbf{f}\colon \Omega\to \mathbb{R}^d$ denotes an external force. The \emph{velocity vector field} 
$\mathbf{v}\colon \overline{\Omega}\to \mathbb{R}^d$ and the scalar \emph{kinematic pressure} $q\colon \Omega\to \mathbb{R}$ are the unknowns in \eqref{eq:fem:main_problem1}. $\mathbf{D}$ denotes the symmetric gradient operator.
Moreover,
we assume that the viscosity of the fluid can be expressed as a function of the \emph{shear rate}. Therefore, we consider a non-degenerate \emph{extra-stress tensor} $\mathbf{S}\colon \mathbb{R}^{d\times d}\to \mathbb{R}_{\textup{sym}}^{d\times d} \, \coloneqq {\{\mathbf{A}\in \mathbb{R}^{d\times d}\mid \mathbf{A}=\mathbf{A}^\top\}}$ that has \emph{$(p,\delta)$-structure} (\textit{cf}.\ Assumption \ref{assum:extra_stress}), where $\delta  >  0$ and $p\in (1,+\infty)$ is the \textit{shear stress exponent}.
A prototypical example for $\mathbf{S}\colon \mathbb{R}^{d\times d}\to \mathbb{R}_{\textup{sym}}^{d\times d}$ is
\begin{align*}
	\mathbf{S}(\mathbf{A}) \coloneqq \nu_0\,(\delta+\vert \mathbf{A}^{\textup{sym}}\vert)^{p-2}\mathbf{A}^{\textup{sym}}\,,
\end{align*}
for $\mathbf{A}\in \mathbb{R}^{d\times d}$, where $\nu_0>0$, $\delta > 0$, $p\in (1,+\infty)$ and $\mathbf{A}^{\textup{sym}}\coloneqq \frac{1}{2}(\mathbf{A}+\mathbf{A}^\top)\in \mathbb{R}^{d\times d}_{\textup{sym}}$.

While the particular case $p=2$ results in the well-known Navier--Stokes equations for \emph{Newtonian} fluids, the cases $p<2$ and $p>2$ correspond to \emph{shear-thinning} and \emph{shear-thickening} fluids, respectively. 
Note that shear-thinning is the most common type of generalized Newtonian behavior of fluids and is seen in many industrial and everyday applications (\textit{cf}.\ \cite[p.\ 61]{bha_dynamics}). This motivates the study of finite element approximations of generalized Navier--Stokes equations \eqref{eq:fem:main_problem1} for small shear stress~exponents~$p<2$.

The system \eqref{eq:fem:main_problem1} is completed by the \textit{inhomogeneous Dirichlet boundary condition}
\begin{alignat}{2}
	\mathbf{v} &= \mathbf{g}_2 &&\quad \textrm{ on } \partial\Omega \,,\label{eq:fem:main_problem20} \\
	\intertext{and, for the sake of generality, the \textit{inhomogeneous divergence constraint}}
	\operatorname{div}\mathbf{v} &= g_1&& \quad \textrm{ in } \Omega\,.\label{eq:fem:main_problem3}
\end{alignat}

This work is dedicated to finite element (FE) approximations of \eqref{eq:fem:main_problem1}--\eqref{eq:fem:main_problem3}, with the focus on convergence rates for velocity and kinematic pressure.\vspace{-2.5mm}

\subsection{Related contributions}
\hspace{5mm}A FE discretization of \eqref{eq:fem:main_problem1} and a proof of weak convergence to a solution of \eqref{eq:fem:main_problem1} is given~in~\cite{DKS13a}. More precisely, in \cite{DKS13a}, homogeneous Dirichlet boundary conditions and the maximal range of the shear stress exponent (\textit{i.e.}, $p\in (\frac{2d}{d+2},+\infty)$)  are considered: in the case $p > \tfrac{2d}{d+1}$, the standard Temam modification (\textit{cf}.\ \cite{Temam84}) of the weak convective term is admissible, while in the case $p > \tfrac{2d}{d+2}$, exactly divergence-free FE pairs can be employed.
A fully-discrete FE discretization of the unsteady Boussinesq system and a proof of weak convergence is given in \cite{FGS22_boussinesq}. There, in the case ${p \in (\tfrac{2d}{d+2}, \tfrac{2d}{d+1}]}$, in order to recover the  crucial cancellation property of the weak convective term,
the usage of a divergence reconstruction operator is proposed. On the theoretical side, this covers and generalizes the case of divergence-free FE pairs. On the \mbox{practical} side, the large computational costs needed by common divergence-free FE pairs are reduced in this way.
To the best of the authors' knowledge, the only contribution on \textit{a priori} error estimates for a finite discretization of \eqref{eq:fem:main_problem1} is \cite{JK23_fem_arxiv}. There, an inhomogeneous Dirichlet boundary condition and the restricted range $p \geq \tfrac{2d}{d+1}$ (for which Temam's modification  (\textit{cf}.\ \cite{Temam84}) is still admissible) is considered.
The inhomogeneous Dirichlet boundary condition is imposed using the Fortin interpolation operator of the velocity FE space, which is (for~many~FE~spaces; we refer the reader to \cite[Rem.\ 2.12]{JK23_fem_arxiv} for a short list) given as a Scott--Zhang interpolation operator plus a divergence correction operator. The divergence correction operator does not need to be computed explicitly as it does not contribute to the Dirichlet boundary traces. Nevertheless,  we still consider the computation of the Scott--Zhang interpolation operator as computationally expensive.\vspace{-2.5mm}

\subsection{New contributions}

\hspace{5mm}The present paper is intended to improve the FE discretization in \cite{JK23_fem_arxiv} in such a way that
the above mentioned gaps are closed. Thus, the new contribution of the present \mbox{paper}~is~\mbox{two-fold}:

\begin{itemize}[noitemsep,topsep=2pt,leftmargin=!,labelwidth=\widthof{$\bullet$}]
	\item  
	First, we vary the discrete Dirichlet boundary condition in such a way that any approximation of the (continuous) Dirichlet boundary condition can be employed. 
	In particular, if the Dirichlet boundary data are sufficiently regular, then this gives a theoretical justification for the usage of nodal interpolation, which is new in this context.	
	
	\item Second, we use the method of divergence-free reconstruction operators in order to treat the full range of the shear stress exponent (\textit{i.e.}, $p\in (\frac{2d}{d+2},+\infty)$).  To the best of the authors' knowledge, this is the first 
	\textit{a priori} error analysis for a discretization of the generalized Navier--Stokes equations \eqref{eq:fem:main_problem1} in the case $p \in (\tfrac{2d}{d+2}, \tfrac{2d}{d+1})$. More precisely, we derive \textit{a priori} error estimates 
	for the velocity vector field with  error decay rates that are optimal for any $p \in (\tfrac{2d}{d+2}, \tfrac{2d}{d+1})$ and  \textit{a priori} error estimates 
	for the kinematic pressure with error decay rates that are optimal for any $p \in (\tfrac{2d}{d+2},2]$. The derivation of \textit{a priori} error estimates 
	for the pressure turns out to be more challenging than in the case of a laminar flow (\textit{cf}. \cite{BBDR12_fem}), \textit{i.e.}, \eqref{eq:fem:main_problem1} without the convective term. More precisely, from the \textit{a priori} error estimates 
	for the velocity vector field and a bootstrap argument, we deduce improved stability estimates for the discrete velocity vector field, which we then employ to prove the \textit{a priori} error estimates for the kinematic pressure.
\end{itemize}
\hspace{5mm}Since these two developments are independent of each other, we decided to present our results in a rather general way. More precisely, we also consider the Temam case, since this --in connection with nodal interpolation of Dirichlet data-- is still new and probably the most relevant case for~application.~Therefore, certain steps of the proofs work out analogously to \cite{JK23_fem_arxiv} and we resort to the respective assertions.\enlargethispage{5mm}

\emph{The paper is structured as follows:} in Section \ref{sec:preliminaries}, the necessary notation and a precise formulation of the continuous problem are given. In Section \ref{sec:discretization}, assumptions on the finite element spaces are stated and the discretization is developed. Section \ref{sec:error_estimates} is dedicated to the proofs of \textit{a priori}  error estimates for the velocity vector field and the kinematic pressure. In Section \ref{sec:fem:num_experiments}, numerical experiments are presented and the theoretical findings are compared to the experimental results.

 \section{Preliminaries}\label{sec:preliminaries}


\subsection{Notation}\label{sub:nota}

\hspace{5mm}Throughout the entire article, we denote by $\Omega \subseteq \mathbb{R}^d$, $d \in \{2,3\}$, either  a bounded polygonal (if $d=2$) or a bounded polyhedral (if $d=3$) Lipschitz domain. We employ $c>0$ to denote a generic constant, that may change from line to line, but does not depend on the crucial quantities. Moreover, we write $u\sim v$ if and only if there exist a constant $c>0$ such that $c^{-1}\, u \leq v\leq c\, u$.

For $k\in \mathbb{N}$ and $p\in [1,\infty]$, we employ the customary Lebesgue spaces ${(L^p(\Omega), \|\cdot\|_{p,\Omega})} $ and Sobolev spaces $(W^{k,p}(\Omega),\|\cdot\|_{k,p,\Omega})$. There exists a linear, continuous, and surjective trace operator $\tr_{\partial\Omega}\colon W^{1,p}(\Omega) \to W^{1-\smash{\frac{1}{p}},p}(\partial\Omega)$. If it gets clear from the context, we do not denote it explicitly for the sake of readability.
The space $\smash{W^{1,p}_{0}(\Omega)}$ is defined as those functions from $W^{1,p}(\Omega)$ whose traces vanish on $\partial\Omega$. 
The conjugate exponent is denoted by $p'  =  \tfrac{p}{p-1}  \in  [1, \infty]$, whereas $p^*$ denotes the critical Sobolev exponent with respect to the embedding $W^{1,p}(\Omega) \hookrightarrow L^{q}(\Omega)$, \textit{i.e.}, if ${p >  d}$, then ${p^*  =  \infty}$, while $p^*$ is used as a placeholder for any finite ${q \in [1, \infty)}$ if $p=d$. If $p=d$, we formally understand $(p^*)'$ as a placeholder for $1+\varepsilon$, where $\varepsilon > 0$ is fixed, but arbitrary. 

We always denote vector-valued functions by lowercase boldface letters and tensor-valued
	functions by capital boldface letters. The Euclidean inner product
between two vectors $\mathbf{a} =(a_1,\dots,a_d),\mathbf{b} =(b_1,\dots,b_d)\in \mathbb{R}^d$ is denoted by 
${\mathbf{a} \cdot\mathbf{b}\coloneqq \sum_{i=1}^d{a_ib_i}}$, while the
Frobenius inner product between two tensors $\mathbf{A}=(A_{ij})_{i,j\in \{1,\dots,d\}},\mathbf{B}=(B_{ij})_{i,j\in \{1,\dots,d\}}\in \mathbb{R}^{d\times d }$ is denoted by
$\mathbf{A}: \mathbf{B}\coloneqq \sum_{i,j=1}^d{A_{ij}B_{ij}}$. For a tensor $\mathbf{A}\in \mathbb{R}^{d\times d}$, we denote its symmetric part by $\mathbf{A}^{\textup{sym}}\coloneqq \frac{1}{2}(\mathbf{A}+\mathbf{A}^\top)\in \mathbb{R}^{d\times d}_{\textup{sym}}\coloneqq {\{\mathbf{A}\in \mathbb{R}^{d\times d}\mid \mathbf{A}=\mathbf{A}^\top\}}$.
Moreover, for a (Lebesgue) measurable set $\omega\subseteq \mathbb{R}^d$  and (Lebesgue) measurable functions $u,v\colon \omega\to \mathbb{R}$ (written $u,v\in  L^0(\omega)$),~we~employ~the~product 
\begin{align*}
	(u,v)_\omega \coloneqq \int_\omega u\, v\,\mathrm{d}x\,,
\end{align*}
whenever the right-hand side is well-defined.
The integral mean of an integrable function $u\in L^1(\omega)$ over a (Lebesgue) measurable set $\omega\subseteq \mathbb{R}^d$ with $\vert \omega\vert>0$ is denoted by 
\begin{align*}
    \langle u\rangle_\omega \coloneqq \frac 1 {|\omega|}\int_\omega u \,\mathrm{d}x\,.
\end{align*} 

For $r \in (1, \infty)$, we abbreviate the function spaces 
\begin{align*}
\begin{aligned}
X^r &\coloneqq  (W^{1,r}(\Omega))^d\,, &&\quad
V^r \coloneqq  (W_{0}^{1,r}(\Omega))^d\,, \\
Y^r &\coloneqq \smash{ L^{r'}(\Omega)}\,, &&\quad
Q^r \coloneqq \smash{ L_0^{r'}(\Omega)} \coloneqq \big \{ z \in \smash{L^{r'}(\Omega) }\mid \langle z \rangle_{\Omega} = 0 \big\}
\,,
\end{aligned}
\end{align*}
as well as  $X\coloneqq X^p$, $V\coloneqq V^p$, $Y\coloneqq Y^p$, and $Q\coloneqq Q^p$.

\subsection{N-functions and Orlicz spaces}

\hspace{5mm}A convex function $\psi \colon \mathbb{R}_{\geq 0} \to \mathbb{R}_{\geq 0}$ is called \emph{N-function} if  $\psi(0)=0$,
$\psi>0$ in $\mathbb{R}_{> 0} $, $\lim_{t\rightarrow0}{\frac{\psi(t)}{t}}=0$ and
${\lim_{t\rightarrow\infty}{\frac{\psi(t)}{t}}=\infty}$. 
The \emph{conjugate N-function} $\psi^*\colon \mathbb{R}_{\geq 0} \to  \mathbb{R}_{\geq 0}$ is defined by
$\psi^*(s)\coloneqq\sup_{t \geq 0}\{s\,t - \psi(t)\}$ for all $s \geq 0$.  
An N-function ${\psi\colon \mathbb{R}_{\geq 0} \to  \mathbb{R}_{\geq 0}}$ satisfies the \emph{$\Delta_2$-condition}, if there exists
$K> 2$ such that ${\psi(2\,t) \leq K\,\psi(t)}$ for all
$t \geq 0$. The smallest such constant is denoted by
$\Delta_2(\psi) > 0$. We need the following refined version of the \textit{$\varepsilon$-Young inequality}: for every $\varepsilon>0$,~there~exists a constant $c_\varepsilon>0 $,
depending only on $\Delta_2(\psi),\Delta_2( \psi ^*)<\infty$, such
that for every ${s,t\geq 0}$, there holds 
\begin{align} \label{eq:fem:orliczyoung}
s\,t&\leq c_\varepsilon \,\psi^*(s)+ \varepsilon \, \psi(t)\,.
\end{align} 

For $p\hspace*{-0.05em} \in \hspace*{-0.05em}(1,\infty)$ and $\delta\hspace*{-0.05em}\geq\hspace*{-0.05em} 0$, we define the \textit{special N-function} $\phi \hspace*{-0.05em}\coloneqq  \hspace*{-0.05em}\phi_{p,\delta}\colon\smash{\mathbb{R}_{\geq 0}\hspace*{-0.05em}\to\hspace*{-0.05em} \mathbb{R}_{\geq 0}}$ by
\begin{align} \label{eq:fem:nfunction} 
\phi(t)\coloneqq  \int _0^t \phi'(s)\, \mathrm{d}s\,, \quad\text{where}\quad
\phi'(t) \coloneqq  (\delta +t)^{p-2} t\,,\quad\textup{ for all }t\geq 0\,.
\end{align}

For an N-function $\psi\colon\mathbb{R}_{\geq 0}\to \mathbb{R}_{\geq 0}$, we define \textit{shifted N-functions} ${\psi_a\colon\mathbb{R}_{\geq 0}\to \mathbb{R}_{\geq 0}}$ via
\begin{align} \label{eq:fem:phi_shifted}
\psi_a(t)\coloneqq  \int _0^t \psi_a'(s)\,  \mathrm{d}s\,,\quad\text{where }\quad
\psi'_a(t)\coloneqq \psi'(a+t)\frac {t}{a+t}\,,\quad\textup{ for all }a, t\geq 0\,.
\end{align}
The shifted special  N-functions $\phi_a\colon \mathbb{R}_{\geq 0}\to \mathbb{R}_{\geq 0}$, $a\ge 0$, and their conjugate N-functions $(\phi_a)^*\colon \mathbb{R}_{\geq 0}\to \mathbb{R}_{\geq 0}$, $a\ge 0$,  satisfy the $\Delta_2$-condition with 
\begin{align}
	\sup_{a\ge 0}{\Delta_2(\phi_a)}&\leq c\,2^{\max \{2,p\}}\,,\\
	\sup_{a\ge 0}{\Delta_2((\phi_a)^*)}&\leq c\,2^{\max \{2,p'\}}\,.
\end{align}
In addition, uniformly in $t,\delta\geq 0$, we have that
\begin{align}
	\phi_a(t) &\sim (\delta + a + t)^{p-2} t^2\,,\label{eq:equ_phi_a}\\
	(\phi_a)^*(t)& \sim ((\delta+a)^{p-1} + t)^{p'-2} t^2\label{eq:equ_phi_a_prime}\,.
\end{align}

Let $\omega\subseteq \mathbb{R}^d$, $d \in \mathbb{N}$, be a (Lebesgue) measurable set. A Carath\'eodory function $\psi \colon \omega \times \mathbb{R}_{\geq 0} \to \mathbb{R}_{\geq 0}$, such that  $\psi(x,\cdot)$ for a.e. $x \in \omega $ is an N-function, is called \emph{generalized N-function}.
For a given generalized N-function $\psi \colon \omega \times \mathbb{R}_{\geq 0} \to \mathbb{R}_{\geq 0}$, the \textit{modular (with respect to $\psi$ on $\omega$)}  of a (Lebesgue) measurable function  $u\in L^0(\omega)$ is defined  by
\begin{align*}
	\rho_{\psi,\omega}(u)\coloneqq \int_\omega
{\psi(\cdot,\vert u\vert)\,\textup{d}x}\,.
\end{align*}

\subsection{The extra-stress tensor and related functions} \hspace{5mm}In this subsection, we specify the assumptions on the extra-stress tensor and important consequences thereof.\enlargethispage{1mm}

\begin{assumption}[extra-stress tensor] \label{assum:extra_stress}
	The extra-stress tensor $\mathbf{S}\colon\mathbb{R}^{d \times d}\to \mathbb{R}^{d \times d}_{\textup{sym}} $ satisfies $ \mathbf{S}\in C^0(\mathbb{R}^{d \times d},\mathbb{R}^{d \times d}_{\textup{sym}} ) $,  $\mathbf{S} (\mathbf{A}) = \mathbf{S} (\mathbf{A}^{\textup{sym}})$ for all ${\mathbf{A}\in \mathbb{R}^{d \times d}}$, and $\mathbf{S}(\mathbf 0)=\mathbf 0$. 
	
Moreover, the extra-stress tensor $\mathbf{S}\colon\mathbb{R}^{d \times d}\to \mathbb{R}^{d \times d}_{\textup{sym}} $ has \textup{$(p,\delta)$-structure}, \textit{i.e.},
	for some $p \in (1, \infty)$, $ \delta\in (0,\infty)$, and the
	N-function $\phi=\phi_{p,\delta}$ defined in \eqref{eq:fem:nfunction}, there exist constants $C_0, C_1 > 0$ such that for every $\mathbf{A},\mathbf{B} \in \mathbb{R}^{d\times d}$, there holds
	\begin{align}
	({\mathbf{S}}(\mathbf{A}) - {\mathbf{S}}(\mathbf{B})) : (\mathbf{A}-\mathbf{B}) &\ge C_0 \,\phi_{\vert \mathbf{A}^{\textup{sym}}\vert}(\vert\mathbf{A}^{\textup{sym}} -
	\mathbf{B}^{\textup{sym}}\vert) \,,\label{assum:extra_stress.1}
	\\
	\vert \mathbf{S}(\mathbf{A}) - \mathbf{S}(\mathbf{B})\vert  &\le C_1 \,
	\phi'_{\abs{\mathbf{A}^{\textup{sym}}}}(\abs{\mathbf{A}^{\textup{sym}} -
		\mathbf{B}^{\textup{sym}}})\,.\label{assum:extra_stress.2}
	\end{align} 
	The constants $C_0,C_1>0$ and $p\in (1,\infty)$ are called the \textup{characteristics} of $\mathbf{S}$.
\end{assumption} 

Since we use estimates that are only available for $\delta>0$ in Section \ref{sec:error_estimates}, we exclude the degenerate case $\delta=0$ in Assumption \ref{assum:extra_stress}.
Closely related to the extra-stress tensor $\mathbf{S}\colon \mathbb{R}^{d \times d} \to \mathbb{R}^{d \times d}_{\textup{sym}} $ with $(p,\delta)$-structure (\textit{cf}.\ Assumption \ref{assum:extra_stress}) is the non-linear mapping $\mathbf{F} 
\colon\mathbb{R}^{d\times d}\to \mathbb{R}^{d\times d}_{\textup{sym}}$, for every $\mathbf{A}\in \mathbb{R}^{d\times d}$  defined by 
\begin{align}
\begin{aligned}
\mathbf{F}(\mathbf{A})&\coloneqq (\delta+\vert \mathbf{A}^{\textup{sym}}\vert)^{\smash{\frac{p-2}{2}}}\mathbf{A}^{\textup{sym}}\,. 
\end{aligned}
\label{eq:def_F}
\end{align}
The non-linear mappings $\mathbf{S},\mathbf{F}\colon \mathbb{R}^{d \times d} \to \mathbb{R}^{d\times d}_{\textup{sym}}$ and
$\phi_a,(\phi_a)^*\colon \mathbb{R}^{\ge 0}\to \mathbb{R}^{\ge 0}$, ${a\ge  0}$,~are~closely~related. 

\begin{proposition}
	\label{lem:growth_SF}
	Let $\mathbf{S}$ satisfy Assumption \ref{assum:extra_stress}, let $\phi$ be defined in \eqref{eq:fem:nfunction} and let $\mathbf{F}$ be defined in \eqref{eq:def_F}. Then, uniformly with respect to 
	$\mathbf{A}, \mathbf{B} \in \mathbb{R}^{d \times d}$, we have that
	\begin{align}\label{eq:growth_S}
	\begin{aligned}
	(\mathbf{S}(\mathbf{A}) - \mathbf{S}(\mathbf{B}))
	:(\mathbf{A}-\mathbf{B} ) &\sim  \abs{ \mathbf{F}(\mathbf{A}) - \mathbf{F}(\mathbf{B})}^2
	\\
	&\sim \phi_{\vert \mathbf{A}^{\textup{sym}}\vert }(\vert \mathbf{A}^{\textup{sym}}
	- \mathbf{B}^{\textup{sym}}\vert )
	\\
	&\sim(\phi_{\vert\mathbf{A}^{\textup{sym}} \vert})^*(\vert\mathbf{S}(\mathbf{A} ) - \mathbf{S}(\mathbf{B} )\vert)\,.
	\end{aligned}
	\end{align}
	The constants in \eqref{eq:growth_S} depend only on the characteristics of ${\mathbf{S}}$.
\end{proposition}
 
 \begin{proof}
 	See \cite[Lem 8.1]{GNSE_DR07}.
 \end{proof}

\subsection{Continuous weak formulation}\vspace{-0.5mm}

\hspace{5mm}In this subsection, we give a precise weak formulation of the  strong formulation \mbox{\eqref{eq:fem:main_problem1}--\eqref{eq:fem:main_problem3}}.\vspace{-0.5mm}

\begin{assumption}\label{assum:fem:data}
	Let $s \coloneqq  s(p) \coloneqq  \max \{p, (\frac{p^*}{2} )'\}$ and let 
	the data satisfy the regularity assumptions $\mathbf{f}\in V^*$, $ g_1  \in L^s(\Omega)$, and $ \mathbf{g}_2\in (W^{\smash{1-\frac{1}{s}},s}(\partial\Omega))^d$ as well as the  compatibility condition\vspace{-0.5mm}
	\begin{align*}
		\int_{\Omega}{g_1 \,\mathrm{d}x} = \int_{\partial\Omega}{\mathbf{g}_2 \cdot\mathbf{n} \, \mathrm{d}s}\,.
	\end{align*} 
\end{assumption}
 
The convective term is consistently reformulated as
\begin{align} \label{eq:fem:conv_reform}
(\operatorname{div}(\mathbf{v} \otimes \mathbf{v}), \mathbf{z})_{\Omega}
= - (\mathbf{v} \otimes \mathbf{v}, \nabla \mathbf{z})_\Omega \,.
\end{align}
Admissibility of the right hand side of \eqref{eq:fem:conv_reform} will be ensured by the restriction $p > \tfrac{2d}{d+2}$ and $\mathbf{z} \in V^s$.

Then, a first weak  formulation of the strong formulation \eqref{eq:fem:main_problem1}--\eqref{eq:fem:main_problem3} is given via:

\textsc{Problem (\hypertarget{Q}{Q}).} Under the  Assumption \ref{assum:fem:data}, find $(\mathbf{v}, q) \in X \times Q^{s}$ with $\mathbf{v} = \mathbf{g}_2$ a.e.\ on $\partial\Omega$ such that for every $(\mathbf{z},z)\in V^{s} \times Y^{s}$, there holds
\begin{align}
(\mathbf{S}(\mathbf{Dv}), \mathbf{Dz})_\Omega - (\mathbf{v} \otimes \mathbf{v}, \nabla \mathbf{z})_\Omega  - (q, \operatorname{div}\mathbf{z})_\Omega  &= (\mathbf{f}, \mathbf{z})_\Omega \,, \label{eq:fem:q_main} \\
(\operatorname{div}\mathbf{v},z)_\Omega &= (g_1,z)_\Omega  \,. \label{eq:fem:q_div} 
\end{align}

We temporarily assume that there exists a divergence and trace lift $\mathbf{g} \in X^s$, which solves the system
\begin{align}\begin{aligned} \label{eq:fem:div_eq}
\operatorname{div}\mathbf{g} &= g_1  &&\quad\text{ in }L^s(\Omega)\,, \\[-0.5mm]
\mathbf{g} &= \mathbf{g}_2  &&\quad\text{ in }(W^{\smash{1-\frac{1}{s}},s}(\partial\Omega))^d\,.
\end{aligned}\end{align}
An explicit construction of such $\mathbf{g}$ is given later in Lemma \ref{lem:fem:gh}. Then, setting $$\mathbf{u} \coloneqq  \mathbf{v} - \mathbf{g} \in V_0\,,$$
where $V_0 \coloneqq  V_0^p$  and for $r\in [1,\infty)$
\begin{align*}
V_0^r \coloneqq \big \{ \mathbf{z} \in V^r \mid\operatorname{div}\mathbf{z} = 0 \text{ a.e.\ in }\Omega\big\}\,,
\end{align*}
we obtain the following weak formulation:

\textsc{Problem (\hypertarget{P}{P}). }
Under the Assumption \ref{assum:fem:data}, find $\mathbf{u} \in V_0$ such that for every $\mathbf{z}\in V_0^{s}$,~there~holds
\begin{align*}
(\mathbf{S}(\mathbf{Du}+\mathbf{Dg}), \mathbf{Dz})_\Omega - ((\mathbf{u}+\mathbf{g}) \otimes (\mathbf{u}+\mathbf{g}), \nabla \mathbf{z})_\Omega
= (\mathbf{f}, \mathbf{z})_\Omega\,.
\end{align*}

First, the existence of the velocity vector field in Problem (\hyperlink{P}{P}) can be guaranteed via the celebrated Lipschitz truncation method (\textit{cf}.\ \cite{dms}). Second, the existence of the kinematic pressure in  Problem (\hyperlink{Q}{Q}) can be inferred on the basis of the following inf-sup stability result. This, in turn, implies equivalence of Problem (\hyperlink{Q}{Q}) and Problem (\hyperlink{P}{P}).\vspace{-0.5mm}

\begin{lemma} \label{lem:fem:cont_infsup}
	Let $r \in (1, \infty)$. Then, for every $z\in Q^r$, there holds
	\begin{align*}
		c\,\|z\|_{r',\Omega}
		\leq  \sup_{\mathbf{z} \in V^r\;:\;\|\mathbf{z}\|_{1,r,\Omega}\leq 1}{(z,\operatorname{div}\mathbf{z})_{\Omega}}\,,
	\end{align*}
	where $c\hspace{-0.1em}>\hspace{-0.1em}0$ depends only on $\Omega$ and $r$. Moreover, $- \nabla \colon \hspace{-0.1em}Q^r \hspace{-0.1em}\to\hspace{-0.1em} (V^r)^*$ is injective and ${\operatorname{img} (-\nabla) \hspace{-0.1em}=\hspace{-0.1em} (\ker \operatorname{div})^{\perp}}$. 
\end{lemma}

\begin{proof}
	See \cite[Thm.\ A.34, Lem. A.40]{EG04_tap}.
\end{proof} 

\section{Discretization} \label{sec:discretization}\enlargethispage{6mm}\vspace{-1mm}

\subsection{Triangulations and finite element spaces}\vspace{-0.5mm}

\begin{assumption}[triangulation] \label{assum:fem:geometry}
	\hspace*{-1mm}We assume that
	$\{\mathcal{T}_h\}_{h>0}$ is a family of conforming
	triangulations of $\overline{\Omega}\subseteq \mathbb{R}^d$,
	$d\in \{2,3\}$, (\textit{cf}.\ \cite{EG21})  consisting of
	\mbox{$d$-dimensional simplices}.
	The parameter $h>0$ refers to the \textup{maximal mesh-size} of $\mathcal{T}_h$, \textit{i.e.}, if $h_K\coloneqq \textup{diam}(K)$ for all $K\in \mathcal{T}_h$, then $h\coloneqq \max_{K\in \mathcal{T}_h}{h_K}$.
	For a simplex $K \in \mathcal{T}_h$,
	we denote the \textup{supremum of diameters of inscribed balls} by $\rho_K>0$ and assume that there exists a \mbox{constant} ${\gamma_0>0}$, independent of $h>0$, such that ${h_K}{\rho_K^{-1}}\le
	\gamma_0$ for all ${K \in \mathcal{T}_h}$. The smallest such constant is called the \textup{chunkiness} of $\{\mathcal{T}_h\}_{h>0}$. For every $K\in \mathcal{T}_h$, the \textit{element patch} is defined as $\omega_K \coloneqq  \{K'\in \mathcal{T}_h\mid K'\cap K \neq \emptyset\}$.
\end{assumption}

\begin{definition}[finite element spaces] \label{defi:fem:spaces}
	The space of polynomials of degree at most $r\in \mathbb{N}\cup\{0\}$ on each simplex is denoted by $\mathbb{P}^r(\mathcal{T}_h)$.
	In addition, for $r\in \mathbb{N}\cup\{0\}$, we set $\mathbb{P}^r_c(\mathcal{T}_h)\coloneqq \mathbb{P}^r(\mathcal{T}_h)\cap C^0(\overline{\Omega})$.
	We assume that the finite element spaces ${X_h \subseteq (\mathbb{P}^m_c(\mathcal{T}_h))^d}$ and $Y_h\hspace*{-0.05em} \subseteq \hspace*{-0.05em}\mathbb{P}^k(\mathcal{T}_h)$ are conforming, \textit{i.e.}, ${X_h\hspace*{-0.05em} \subseteq \hspace*{-0.05em} X}$ and ${Y_h\hspace*{-0.05em} \subseteq\hspace*{-0.05em} Y}$. Then, we define the spaces\vspace*{-0.5mm}
	\begin{align*}
		V_h &\coloneqq  X_h \cap V\,,\\
		Q_h& \coloneqq  Y_h \cap Q\,,\\ 
		V_{h,0}&\coloneqq \big\{ \mathbf{z}_h \in V_h \mid (\mathrm{div}\,\mathbf{z}_h,y_h)_\Omega = 0 \textrm{ for all } y_h \in Q_h \big\}\,.
	\end{align*} 
\end{definition}

The following assumption is needed in the construction of a discrete (divergence-corrected) Lipschitz truncation  (\textit{cf}.\ \cite{DKS13a} or \cite[Ass.\ 2.20]{Tscherpel_phd}).
\begin{assumption}[locally supported basis of $Y_h$] \label{assum:fem:localbasis}
	We assume that $Y_h$ has a locally supported basis $Y_h = \textup{span} \{q_h^1, \dots, q_h^l \}$ such that $q^i_h|_K \neq 0$ implies  that $\textup{supp}\, q_h^i \subseteq \omega_K$ for all $i \in \{1, \dots, l\}$ and $K \in \mathcal{T}_h$.
\end{assumption}

\begin{assumption}[projection for $Y_h$] \label{assum:fem:projection_y}
	We suppose that  $\mathbb{R}\subseteq Y_h$ 
	and that there exists a linear projection operator $\Pi_h^Y \colon Y \to Y_h$, \textit{i.e.}, $\Pi_h^Yz_h=z_h$ for all $z_h\in Y_h$, such that for every $z \in Y$ and ${K \in \mathcal{T}_h}$, there holds\vspace*{-0.5mm}
	\begin{align*}
		\langle \vert \Pi_h^Y\! z\vert \rangle_K\leq c\, \langle \vert z\vert \rangle_{\omega_K}\,.
	\end{align*}
\end{assumption}

\begin{assumption}[projection for $X_h$] \label{assum:fem:projection_x}
	We suppose that $(\mathbb{P}^1_c(\mathcal{T}_h))^d \subseteq X_h$ and that there exists
	a linear projection operator $\Pi_h^X \colon X \to X_h$, \textit{i.e.}, $\Pi_h^X\mathbf{z}_h=\mathbf{z}_h$ for all $\mathbf{z}_h\in X_h$,~with~the~following~properties:
	\begin{itemize}[noitemsep,topsep=2pt,leftmargin=!,labelwidth=\widthof{(iii)}]
		\item[(i)] \textup{Local $W^{1,1}$-stability:}  \textcolor{black}{For} every $\mathbf{z} \in X$ and $K \in \mathcal{T}_h$, there holds 
		\begin{align*}
			\langle \vert\Pi_h^X \mathbf{z}\vert \rangle_K  \leq c\,\langle \vert\mathbf{z}\vert \rangle_{\omega_K} + c\, h_K \,\langle \vert \nabla \mathbf{z}\vert  \rangle_{\omega_K} \,.
		\end{align*}
		\item[(ii)] \textup{Preservation of zero boundary values:} there holds $\Pi_h^X(V) \subseteq V_h$,
		\item[(iii)] \textup{Preservation of divergence in the $Y_h^*$-sense:} For every $\mathbf{z} \in V$ and  $z_h \in Y_h$, there holds 
		\begin{align*}
			(\mathrm{div}\,\mathbf{z}, z_h)_\Omega = (\mathrm{div}\,\Pi_h^X \mathbf{z}, z_h)_\Omega \,.
		\end{align*}
	\end{itemize}
\end{assumption}

\begin{remark}
	A projection operator $\Pi_h^X\colon X\to X_h$ which satisfies Assumption \ref{assum:fem:projection_x} is traditionally referred to as Fortin interpolation operator. Its existence for typical finite element pairs is discussed, \textit{e.g.}, in \cite{ET} (see also \cite[p.~23]{GazcaGmeinederMaringovaTscherpel2025}).
\end{remark}

Assumption \ref{assum:fem:projection_x} implies the following discrete inf-sup stability result.
\begin{lemma}[discrete inf-sup condition] \label{lem:fem:discrete-inf-sup}
	Let Assumption \ref{assum:fem:projection_x} be satisfied and $r \in (1, \infty)$. Then, for every $z_h \in Q_h$,
	there holds\vspace*{-0.5mm}
	\begin{align*}
		c\,\|z_h\|_{r', \Omega} &\leq c \sup_{\mathbf{z}_h  \in V_h\;:\;\|\mathbf{z}_h\|_{1,r,\Omega}\leq 1}{(z_h,\mathrm{div}\,\mathbf{z}_h)_{\Omega}}\,, 
	\end{align*}
	where the constant $c>0$ depends on $r$, $\gamma_0$, $\Omega$ and the choice of finite element spaces.
\end{lemma}

\begin{proof}
	See \cite[Lem.\ 4.1]{BBDR12_fem}.
\end{proof}

\subsection{Discrete Dirichlet boundary condition}

\hspace*{5mm}We assume that among the degrees of freedom  parametrizing $X_h$, those which are relevant for the Dirichlet boundary \mbox{condition} \eqref{eq:fem:main_problem20} can be uniquely determined (\textit{cf}. \cite[Sec. 19.4]{EG21}). We interpret the Dirichlet boundary condition \eqref{eq:fem:main_problem20} as Dirichlet boundary condition at the Dirichlet degrees of \mbox{freedom} in $X_h$. As a consequence, an interpolation of boundary data $\mathbf{g}_b \in X^s$ needs to be computed. The Fortin interpolation operator has ideal analytical properties for that purpose (\textit{cf}.\   \cite{JK23_fem_arxiv} and Remark \ref{rem:fem:classical_interpolation}), but this comes with significant computational costs. Therefore, we develop a setting which allows for interpolation operators which are computationally inexpensive and straightforward to implement. Since regularity of boundary data is rarely critical in applications, the nodal interpolation operator is probably the most interesting example.
In this subsection, we formulate requirements on continuous and discrete boundary data. Then, we discuss the applicability of some interpolation operators.\newpage

The following assumptions will be sufficient for existence and convergence of discrete solutions:
\begin{assumption}\label{assum:fem:projection_x2_weak}
	We assume that $\mathbf{g}_b \in X^s$ satisfies $\mathbf{g}_b = \mathbf{g}_2$ a.e. on $\partial\Omega$ and that $\mathbf{g}_b^h \in X_h$ satisfy
	\begin{align*}
		\mathbf{g}_b^h \to \mathbf{g}_b \quad \textrm{in } X^s \quad (h \to 0)\,.
	\end{align*}
\end{assumption}

The regularity $\mathbf{g}_b \in X^s$ complies with the best known assumptions for existence of weak solutions to Problem (\hyperlink{Q}{Q}) (at least for small values of $p$), \textit{cf}.\ \cite{JR21_inhom}.

Higher regularity of $\mathbf{g}_b\in X^s$ and approximability by $\mathbf{g}_b^h\in X_h$ will be employed in the proof of \textit{a priori} error estimates in Section \ref{sec:error_estimates}:
\begin{assumption}\label{assum:fem:projection_x2}
	We assume that there exist $k \geq 2$, $\beta \in [1, \infty)$, $\mathbf{g}_b \in (W^{k,\beta}(\Omega))^d$, and $\mathbf{g}_b^h \in X_h$ such that the following conditions are satisfied:
	\begin{itemize}[noitemsep,topsep=2pt,leftmargin=!,labelwidth=\widthof{(iii)}]
		\item[(i)] $(W^{k,\beta}(\Omega))^d \hookrightarrow  X^s$ and $(\mathbb{P}^{k-1}_c(\mathcal{T}_h))^d \subseteq X_h$;
		\item[(ii)] $\mathbf{g}_b = \mathbf{g}_2$ a.e.\ on $\partial\Omega$;
		\item[(iii)] 
		$\mathbf{g}_b \in (W^{k,\smash{\smash{\widetilde{\beta}}}}(\Omega))^d$ for some $\smash{\smash{\widetilde{\beta}}} \in [\beta, \infty)$ implies that 
		\begin{align} \label{eq:fem:proj_x2_sappr}
			\|\mathbf{g}_b^h - \mathbf{g}_b\|_{1,\smash{\smash{\widetilde{\beta}}},\Omega}
			\leq c \,h^{k-1} \|\nabla^k\mathbf{g}_b\|_{\smash{\smash{\widetilde{\beta}}}, \Omega}\,.
		\end{align}
	\end{itemize}
\end{assumption}

Note that Assumption \ref{assum:fem:projection_x2} implies Assumption \ref{assum:fem:projection_x2_weak} (\textit{cf}.\  reasoning in \cite[p. 35]{EG21-2}). While $k=2$, $\beta=p$ in Assumption \ref{assum:fem:projection_x2} is certainly a natural choice in order to establish linear convergence rates, different regularities may be imposed due to the admissibility of interpolation operators and technical restrictions in the analysis.

\begin{lemma}[classical interpolation operators] \label{lem:fem:proj_x2_w11}
	Under Assumption \ref{assum:fem:projection_x2}(i), let $\tilde{\Pi}_h^X\colon X \to  (\mathbb{P}^{k-1}_c(\mathcal{T}_h))^d$ be a linear projection operator 
	that is \emph{quasi-locally \mbox{$W^{1,1}$-stable}}, \textit{i.e.}, for every $\mathbf{z} \in (W^{1,1}(\Omega))^d$ and $K\in \mathcal{T}_h$, there holds
	\begin{align} \label{eq:fem:projection_x2_assum}
		\langle \vert\tilde{\Pi}_h^X \mathbf{z}\vert \rangle_K  \leq c\,\langle \vert\mathbf{z}\vert \rangle_{\omega_K} + c\, h_K\,\langle  \vert \nabla \mathbf{z}\vert  \rangle_{\omega_K}\,.
	\end{align}
	Then, Assumption \ref{assum:fem:projection_x2}(iii) is satisfied
	for  $\mathbf{g}_b^h \coloneqq \tilde{\Pi}_h^X\mathbf{g}_b\in  X_h$ with $\smash{\smash{\widetilde{\beta}}}=\beta$. 
\end{lemma}
\begin{proof} 
	See \cite[Folg.\ 3.28]{Dziuk_tun}.  
\end{proof}

\begin{remark}[classical interpolation operators] \label{rem:fem:classical_interpolation}
	Suppose that Assumption \ref{assum:fem:projection_x2}(i) is satisfied. Then, 
	the Fortin interpolant 
	(\textit{cf}.\   Assumption \ref{assum:fem:projection_x}), the Clément \cite{Clement75_interpolation}, and the Scott-Zhang interpolant \cite{SZ90_interpolation} or the generic interpolant from \cite{EG21}, due to Lemma \ref{lem:fem:proj_x2_w11}, satisfy Assumption \ref{assum:fem:projection_x2}(iii).
\end{remark}

\begin{remark}[(global) $L^2$-projection operator]
	Under Assumption \ref{assum:fem:projection_x2}(i) and if
	the family of triangulations $\{\mathcal{T}_h\}_{h>0}$ is quasi-uniform or appropriately graded, the (global) $L^2$-projection operator $\Pi_{h,L^2}^X\colon  X\to  (\mathbb{P}^{k-1}_c(\mathcal{T}_h))^d$ is (globally) $L^{1}(\Omega)$-stable \cite{Douglas74, DST21_L2} and, hence, satisfies Assumption~\ref{assum:fem:projection_x2}(iii) \cite[Folg.\ 3.28]{Dziuk_tun}.
\end{remark}

The next Lemma shows that despite the nodal interpolation operator $\Pi_{h, N}^X$ is not $W^{1,1}(\Omega)$-stable, it still satisfies Assumption \ref{assum:fem:projection_x2} (under appropriate regularity assumptions).

\begin{lemma}[nodal interpolation] \label{lem:fem:proj_x2_nodal}
	Assume that $k \in \mathbb{N}$ with $k\ge 2$ and ${\beta \in [1, \infty)}$ are such that $W^{k,\beta}(\Omega) \hookrightarrow  C^0(\overline{\Omega})$ and $(\mathbb{P}^{k-1}_c(\mathcal{T}_h))^d \subseteq X_h$. We suppose that each degree of freedom of $X_h$ is represented by evaluation at a certain point. Let $\Pi_{h, N}^X\colon (W^{k,\beta}(\Omega))^d\to (\mathbb{P}^{k-1}_c(\mathcal{T}_h))^d$ be the nodal interpolant.
	Then, 
	\eqref{eq:fem:proj_x2_sappr} is satisfied for $\mathbf{g}_b^h \coloneqq  \Pi_{h,N}^X \mathbf{g}_b\in X_h $.
\end{lemma}
\begin{proof}
	See \cite[Satz 3.31]{Dziuk_tun}. 
\end{proof}

\begin{lemma}\label{lem:fem:projection_x2_oappr}
	Let Assumption \ref{assum:fem:projection_x2} be satisfied. Moreover, let $\delta>0$, $\mathbf{z} \in X$, and $\mathbf{g}_b \in W^{k,\max\{2,p,\beta\}}(\Omega)$. Then, there holds
	\begin{align*}
		\rho_{\phi_{\vert\mathbf{Dz}\vert} ,\Omega} (\vert \mathbf{Dg}_b^h - \mathbf{Dg}_b\vert) 
		\leq c\, h^2\, \|\nabla^k \mathbf{g}_b\|_{\max\{2,p,\beta\}, \Omega}^2\,,
	\end{align*}
	where $c>0$ depends on $\delta^{-1}$, $\|\mathbf{g}_b\|_{1,p,\Omega}$, and $\|\mathbf{z}\|_{1,p,\Omega}$. In addition,  
	the assertion remains true if $\mathbf{g}_b^h$ is replaced by $\Pi_h^X \mathbf{g}_b$.
\end{lemma}
\begin{proof}
	If  $p<2$ and $\delta > 0$, the estimate follows from $\phi_{\vert \mathbf{Dz}\vert}(t) \sim (\delta + \vert\mathbf{Dz}\vert + t)^{p-2} t^2$ $ \leq \delta^{p-2} t^2$ uniformly in $t\ge 0$ and Sobolev approximability (\textit{cf}. Assumption \ref{assum:fem:projection_x2}(iii)).
	If $p \geq 2$, we use that $\phi_{\vert\mathbf{Dz}\vert}(t) \sim (\delta + \vert\mathbf{Dz}\vert + t)^{p-2} t^2$ uniformly in $t\ge 0$, Hölder's inequality, $\mathbf{g}_b^h \to \mathbf{g}_b$ in $X$ $(h\to 0)$, and Assumption \ref{assum:fem:projection_x2}(iii) to estimate\vspace{-0.5mm}
	\begin{align*}
		\rho_{\phi_{\vert\mathbf{Dz}\vert} ,\Omega} (\vert \mathbf{Dg}_b^h - \mathbf{Dg}_b\vert) 
		&\leq c\,\|(\max \{\delta, \vert\mathbf{Dz}\vert, \vert \mathbf{Dg}_b^h - \mathbf{Dg}_b\vert \})^{\smash{p-2}}\vert\mathbf{Dg}_b^h - \mathbf{Dg}_b\vert^2 \|_{1,\Omega} \\
		&\leq c\,\|\max \{ \delta, \vert\mathbf{Dz}\vert, \vert\mathbf{Dg}_b^h - \mathbf{Dg}_b\vert \}\|_{p,\Omega}^{p-2} \|\mathbf{Dg}_b^h - \mathbf{Dg}_b\|_{p,\Omega}^2 \\
		&\leq c \, h^2 \|\nabla^k \mathbf{g}_b\|_{\max\{p,\beta\},\Omega}^2 \,.
	\end{align*}
	The assertion for $\mathbf{g}_b^h$ replaced by $\Pi_h^X \mathbf{g}_b$ follows in virtue of Remark \ref{rem:fem:classical_interpolation}.
\end{proof}

We note that the difference $\Pi_h^X \mathbf{g}_b - \mathbf{g}_b^h$ approaches zero with an at least linear rate in $h$: due to Hölder's inequality, Sobolev approximability (\textit{cf}.\ \cite[Thm. 4.6]{DR07_interpolation}), and \eqref{eq:fem:proj_x2_sappr}, for every $t \in [1, \infty)$ and $\smash{\widetilde{\beta}} \coloneqq \max \{t, \beta\}$, there holds\vspace{-0.5mm}
\begin{align}\begin{aligned} \label{eq:fem:ph-gb}
		\|\Pi_h^X \mathbf{g}_b - \mathbf{g}_b^h\|_{1,t,\Omega}
		&\leq c \, \|\Pi_h^X \mathbf{g}_b - \mathbf{g}_b\|_{1,\smash{\widetilde{\beta}},\Omega} + c \, \|\mathbf{g}_b - \mathbf{g}_b^h\|_{1,\smash{\widetilde{\beta}},\Omega}
		\\&\leq c\, h\, \|\nabla^k \mathbf{g}_b\|_{\smash{\widetilde{\beta}},\Omega} \,.
\end{aligned}\end{align}

\subsection{Discrete convective term}\enlargethispage{6.5mm}\vspace{-0.5mm}

\hspace{5mm}Classically, skew-symmetry plays a crucial role in the analysis of the convective term. While this property seems out of reach in the case of a non-zero divergence constraint on the velocity, we may at least preserve certain cancellation properties. Therefore, in the case $p \geq \smash{\tfrac{2d}{d+1}}$, we employ discretizations ${b,\widetilde{b}\colon (X)^2\times X^s \to\mathbb{R}}$, commonly referred to as \textit{Temam's modifications}, for every ${(\mathbf{u},\mathbf{v},\mathbf{w})\in (X)^2\times X^s}$, (\textit{cf}.\ \cite{Temam84}) defined by 
\begin{align}\label{eq:fem:defi_conv_temam}
	\begin{aligned}
		b(\mathbf{u},\mathbf{v},\mathbf{w}) 
		&\coloneqq  \tilde{b}(\mathbf{u},\mathbf{v},\mathbf{w}) + \tfrac{1}{2} (g_1 \mathbf{u}, \mathbf{w})_\Omega \\
		&\coloneqq  \tfrac{1}{2} (\mathbf{w}\otimes\mathbf{u}, \nabla \mathbf{v} + g_1 \mathbf{I})_\Omega - \tfrac{1}{2} (\mathbf{v}\otimes\mathbf{u}, \nabla \mathbf{w})_\Omega\,.
	\end{aligned}
\end{align}

\begin{lemma}\label{lem:fem:admissible}
	Let $p \geq \smash{\tfrac{2d}{d+1}}$. Then, Temam's modifications ${b,\widetilde{b}\colon (X)^2\times X^s \to \mathbb{R}}$ are well-defined and bounded. Moreover, $b\colon (X)^2\times X^s \to \mathbb{R}$
	is consistent, \textit{i.e.}, for every $\mathbf{v} \in \textcolor{black}{X}$ with $\mathrm{div}\,\mathbf{v} = g_1$ a.e.\ in $\Omega$ and ${\mathbf{z}\in \textcolor{black}{V^s}}$, there holds
	\begin{align} \label{eq:fem:temam_consistent}
		b(\mathbf{v}, \mathbf{v}, \mathbf{z}) = - (\mathbf{v}\otimes\mathbf{v}, \nabla\mathbf{z})_{\Omega}\,.
	\end{align}
	For every $\mathbf{v} \in X$ and $\mathbf{z} \in V^s$, there holds\vspace{-0.5mm}
	\begin{align} \label{eq:fem:temam_skew} 
			\tilde{b}(\mathbf{v}, \mathbf{z}, \mathbf{z}) = 0\,.
	\end{align} 
\end{lemma}
\begin{proof}
	See \cite[Lem.\ 2.16]{JK23_fem_arxiv}.
\end{proof}

The lower bound on the range of $p$ is clearly suboptimal, as the (continuous) existence theory~for weak solutions allows for values $p > \tfrac{2d}{d+2}$ (\textit{cf}.\ \cite{dms}). The underlying issue can be circumvented~by~the~usage~of~a \emph{(divergence) reconstruction operator}, which maps discretely 
to~exactly~\mbox{divergence-free}~functions.\vspace{-0.5mm}

\begin{definition}[reconstruction operator] \label{defi:fem:reconstruction}
	Let $Z_h \subset H(\operatorname{div};\Omega)$\footnote{$H(\operatorname{div};\Omega)\coloneqq\{\mathbf{z}\in (L^2(\Omega))^d\mid \operatorname{div}\mathbf{z}\in L^2(\Omega)\}$.} be a finite element space  such that $Z_h|_K \subseteq (W^{1,\infty}(K))^d$ for all $K \in \mathcal{T}_h$.
	A linear operator $\Sigma_h\colon X \to X_h + Z_h$ is called a \emph{(divergence) reconstruction operator}, if the following conditions are satisfied: 
	\begin{itemize}[noitemsep,topsep=2pt,leftmargin=!,labelwidth=\widthof{(iii)}]
		\item[(i)] \emph{(divergence preservation)} For  every $\mathbf{z} \in X$ with $(\mathrm{div}\,\mathbf{z}, y_h)_{\Omega} = 0$ for all ${y_h \in Y_h}$, there holds
		\begin{align*}
			\mathrm{div}\,\Sigma_h \mathbf{z} = 0\quad \text{ a.e.\ in }\Omega\,.
		\end{align*}
		\item[(ii)] \emph{(consistency)} For every $m \in \{0,1\}$, $r \in [1, \infty]$, $K \in \mathcal{T}_h$, and $\mathbf{z} \in X^r$, there holds 
		\begin{align*} 
			\norm{\mathbf{z} - \Sigma_h \mathbf{z}}_{m,r,K}
			\leq c\, \smash{h_K^{1-m}}\, \|\nabla\mathbf{z}\|_{r,K} \,.
		\end{align*}
		\item[(iii)] \emph{(stability of divergence)} For every  $\mathbf{z} \in X$ with $\mathrm{div}\,\mathbf{z} \in L^r(\Omega)$, where $r \in [1, \infty)$, there holds 
		\begin{align*}
			\norm{\mathrm{div}\,\Sigma_h \mathbf{z}}_{r,\Omega} \leq c \norm{\mathrm{div}\,\mathbf{z}}_{r,\Omega} \,.
		\end{align*}
	\end{itemize}
\end{definition}

\begin{lemma} \label{lem:fem:sigma_stability}
	Let $r \in [1, \infty]$. If  $\Sigma_h\colon X \to X_h + Z_h$  is a (divergence) reconstruction operator in the sense of 
	Definition \ref{defi:fem:reconstruction}, then, for every $\mathbf{z} \in X^r$, there holds
	\begin{align*}
		\|\Sigma_h \mathbf{z}\|_{r^*,\Omega} \leq c\, \|\mathbf{z}\|_{1,r,\Omega}\,.
	\end{align*}
\end{lemma}
\begin{proof}
	Due to the triangle inequality and a Sobolev embedding in $\Omega$, there holds
	\begin{align}\label{lem:fem:sigma_stability.1}
		\begin{aligned} 
		\|\Sigma_h \mathbf{z}\|_{r^*,\Omega}&\leq c\, \|\mathbf{z}\|_{r^*,\Omega} + \|\mathbf{z}-\Sigma_h \mathbf{z}\|_{r^*,\Omega}
		\\&\leq c\, \|\mathbf{z}\|_{1,r,\Omega} + \|\mathbf{z}-\Sigma_h \mathbf{z}\|_{r^*,\Omega}\,.
	\end{aligned}
	\end{align}
	Let $K \in \mathcal{T}_h$ be fixed, but arbitrary. A Sobolev embedding in the unit simplex and a transformation to $K$ (\textit{cf}.\ \cite[Folg.\ 3.30]{Dziuk_tun}) yields the local Sobolev inequality
	\begin{align}\label{lem:fem:sigma_stability.2}
		\|\mathbf{z}-\Sigma_h \mathbf{z}\|_{r^*,K} \leq c\, h_K^{-1} \|\mathbf{z}-\Sigma_h \mathbf{z}\|_{r,K} + c\,\|\nabla\mathbf{z}-\nabla\Sigma_h \mathbf{z}\|_{r,K}\,.
	\end{align}
	 With consistency (\textit{cf}.\ Definition \ref{defi:fem:reconstruction}(ii)), from \eqref{lem:fem:sigma_stability.2}, it follows that
	\begin{align}\label{lem:fem:sigma_stability.3}
		\|\mathbf{z}-\Sigma_h \mathbf{z}\|_{r^*,K}
		\leq c \,\|\mathbf{z}\|_{1,r,K}\,.
	\end{align}

	If $r^* < \infty$, we deduce the assertion from \eqref{lem:fem:sigma_stability.1} by summation of \eqref{lem:fem:sigma_stability.3} together with the power mean inequality (\textit{cf}.\ \cite[Ex.\ 12.1]{EG21}).  
	
	Else if $r^* = \infty$, we conclude by \eqref{lem:fem:sigma_stability.1}, \eqref{lem:fem:sigma_stability.3} and $\| \mathbf{y}\|_{r^*,\Omega} \leq \max_{K \in \mathcal{T}_h} \| \mathbf{y}\|_{r^*,K}$~for~all~${\mathbf{y} \in X^r}$.
\end{proof}

Next, we present short lists of common mixed finite element spaces $\{X_h\}_{h>0}$ and $\{Y_h\}_{h>0}$ with projectors $\{\Pi_h^X\}_{h>0}$ and $\{\Pi_h^Y\}_{h>0}$ on regular triangulations $\{\mathcal{T}_h\}_{h>0}$ satisfying both  Assumption \ref{assum:fem:projection_y} and Assumption \ref{assum:fem:projection_x}, for which
reconstruction operators $\{\Sigma_h\}_{h>0}$ in the sense of Definition~\ref{defi:fem:reconstruction}~are~available.

\begin{remark}[discontinuous pressure space]\label{rem:fem:reconstruction}
	We distinguish between element-wise constant (\textit{i.e.}, $k=0$) and element-wise  affine (\textit{i.e.}, $k=1$)  pressure spaces ${Y_h=\smash{\mathbb{P}^k(\mathcal{T}_h)}}$: 
	\begin{itemize}[noitemsep,topsep=2pt,leftmargin=!,labelwidth=\widthof{(ii)}]
		\item[(i)] \emph{Element-wise constant pressure space:} ($k=0$)
		\begin{itemize}[noitemsep,topsep=2pt,leftmargin=!,labelwidth=\widthof{(i.a)}]
			\item[(i.a)]  The \emph{first-order Bernardi--Raugel  element} (\textit{cf}.\ \cite{BR85}) for $d\in \{2,3\}$, \textit{i.e.}, 
			$X_h =(\mathbb{P}^1_c(\mathcal{T}_h)\bigoplus \mathbb{B}_{\tiny \mathscr{F}}(\mathcal{T}_h))^d$, where $\mathbb{B}_{\tiny \mathscr{F}}(\mathcal{T}_h)$ is the facet bubble \mbox{function} space. 
			
			\item[(i.b)] The \emph{$\mathbb{P}^2$-$\mathbb{P}^0$-element} for $d=2$, \textit{i.e.}, $X_h=(\mathbb{P}^2_c(\mathcal{T}_h))^2$.
		\end{itemize} 
		\item[(ii)] \emph{Element-wise linear pressure space:} ($k=1$)
		\begin{itemize}[noitemsep,leftmargin=!,labelwidth=\widthof{(ii.a)}]
			
			\item[(ii.a)]  The \emph{conforming Crouzeix--Raviart element} (\textit{cf}.\ \cite{CR73}) for $d=2$, \textit{i.e.}, $X_h=(\mathbb{P}^2_c(\mathcal{T}_h)\bigoplus\mathbb{B}(\mathcal{T}_h))^2$, where $\mathbb{B}(\mathcal{T}_h)$ is the bubble function space.

			\item[(ii.b)] The  \emph{second-order  Bernardi--Raugel  element} (\textit{cf}.\ \cite{BR85}) for  $d=3$, \textit{i.e.}, $X_h=(\mathbb{P}^2_c(\mathcal{T}_h)\bigoplus \mathbb{B}_{\tiny \mathscr{F}}(\mathcal{T}_h)\linebreak\bigoplus\mathbb{B}(\mathcal{T}_h))^3$.
		\end{itemize}
	\end{itemize} 
	For $k\in \{0,1\}$, the reconstruction, \textit{e.g.}, can be done in the Raviart--Thomas \mbox{element} $Z_h \coloneqq  \mathcal{R}T^k(\mathcal{T}_h)$ with $\Sigma_h\coloneqq \Pi_h^{rt,k}\colon H(\operatorname{div};\Omega) \cap (L^{2+\varepsilon}(\Omega))^d \to  \mathcal{R}T^k(\mathcal{T}_h)$, $\varepsilon>0$, given via the Fortin interpolation operator, for every $\mathbf{z}\in H(\operatorname{div};\Omega) \cap (L^{2+\varepsilon}(\Omega))^d$ \mbox{defined} by\vspace*{-1mm}
	\begin{align}
		\label{eq:fortin_rt.1}
		\int_F{\Pi_h^{rt,k} \mathbf{z}\cdot \mathbf{n}\,\varphi_{F} \,\mathrm{d}s} &=	\int_F{ \mathbf{z}\cdot \mathbf{n}\, \varphi_F\,\mathrm{d}s}&&\text{for all }\varphi_F\in \mathbb{P}^k(F)&&\text{for all }F\in \mathcal{F}_h\,, \\
		\label{eq:fortin_rt.2}
		\int_K{\Pi_h^{rt,k} \mathbf{z}\cdot \mathbf{\varphi}_K \,\mathrm{d}x} &=	\int_K{ \mathbf{z}\cdot \mathbf{\varphi}_K\,\mathrm{d}x}&&\text{for all }\mathbf{\varphi}_K\in (\mathbb{P}^{k-1}(K))^d&&\text{for all }K\in \mathcal{T}_h\,,\\[-6mm]\notag
	\end{align}
	where $\mathcal{F}_h\hspace{-0.1em}\coloneqq\hspace{-0.1em} \{K\cap K'\mid K,K'\hspace{-0.1em}\in\hspace{-0.1em} \mathcal{T}_h,\,\textup{dim}_{\mathscr{H}}(K\cap K')=d-1\}$\footnote{For a set $M\subseteq \mathbb{R}^d$, by $\textup{dim}_{\mathscr{H}}(M)\coloneqq \inf\{s>0\mid \mathscr{H}^s(M)=0\}$ we denote the \emph{Hausdorff dimension}.} is the set of edges (if ${d\hspace{-0.1em}=\hspace{-0.1em}2}$)~or~facets~(if~${d\hspace{-0.1em}=\hspace{-0.1em}3}$). In particular, in relation \eqref{eq:fortin_rt.2}, we use the convention $\mathbb{P}^{-1}(K)\coloneqq \emptyset$, such that  in the case $k=0$, relation \eqref{eq:fortin_rt.2} is trivially satisfied.  Moreover, there holds $\mathrm{div}\,Z_{h}  = Y_{h}$ (this structural connection is essential for the method), which then implies that discretely divergence-free functions are mapped to exactly divergence-free functions (\textit{cf}. \cite[Eq. (5.8)]{JLMNR17_on-the-divergence} for the proof in the second case; the first case works analogously). The stability and approximation properties in Definition \ref{defi:fem:reconstruction}~are~proved~in~\mbox{\cite[Thm.~16.4]{EG21}}.
\end{remark}

\begin{remark}[continuous pressure space]
	We restrict to an element-wise affine, globally continuous pressure space $\smash{Y_h=\mathbb{P}^1_c(\mathcal{T}_h)}$: 
	\begin{itemize}[noitemsep,topsep=2pt,leftmargin=!,labelwidth=\widthof{(ii)}]
		\item[(i)]    The \emph{MINI element} (\textit{cf}.\ \cite{ABF84}) for $d\in \{2,3\}$, \textit{i.e.}, $X_h=(\mathbb{P}^1_c(\mathcal{T}_h)\bigoplus\mathbb{B}(\mathcal{T}_h))^d$.
		
		\item[(ii)] The \emph{Taylor--Hood element} (\textit{cf}.\ \cite{TH73}) for $d\in\{2,3\}$, \textit{i.e.}, $X_h=(\mathbb{P}^2_c(\mathcal{T}_h))^d$.
	\end{itemize} 
	For these two examples, 
	a more involved construction yields reconstruction operators (\textit{cf}.\ \cite{LLMS17_divergence-free}). The proof of the properties in Definition \ref{defi:fem:reconstruction} is similar to \cite[Thm.\ 12]{LLMS17_divergence-free}. However, it remains an open question whether they fulfill the stronger assumptions in Definition~\ref{defi:fem:reconstruction}.
\end{remark}


\begin{assumption}\label{assum:fem:reconstruction}
	We assume that $p \geq \smash{\tfrac{2d}{d+1}}$ or that $p  > \smash{\tfrac{2d}{d+2}}$ together with the existence of a reconstruction operator $\Sigma_h\colon X\to X_h+Z_h$ in the sense of Definition \ref{defi:fem:reconstruction}.
\end{assumption}

In the second case in Assumption \ref{assum:fem:reconstruction},  the \emph{discrete convective term}, for every $\smash{(\mathbf{u},\mathbf{v},\mathbf{w})}\in 
(X)^2\times X$, is defined by 
\begin{align}\label{eq:fem:defi_conv_recon}
		b(\mathbf{u},\mathbf{v},\mathbf{w})  
		\coloneqq  - (\mathbf{v}\otimes \Sigma_h \mathbf{u}, \nabla \mathbf{w})_{\Omega}\,. 
\end{align}

\begin{lemma}\label{lem:fem:admissible2}
	Let the second case in Assumption \ref{assum:fem:reconstruction} be satisfied. Then, the discrete convective term $b\colon (X)^2\times X\to \mathbb{R}$ is well-defined, bounded, 
	 approximately consistent, i.e., for every $\mathbf{v} \in X$ and $\mathbf{z} \in V^s$, we have that\vspace*{-0.75mm}
	\begin{align}  \label{eq:fem:conv-recon-consistency}
			b(\mathbf{v}, \mathbf{v}, \mathbf{z}) 
			= - (\mathbf{v}\otimes\mathbf{v}, \nabla\mathbf{z})_{\Omega}
			+ (\mathbf{v}\otimes \{\mathbf{v} - \Sigma_h \mathbf{v}\}, \nabla\mathbf{z})_{\Omega} \,,\\[-5.5mm]\notag
		\end{align}
	and  has the cancellation property, i.e., for every $\mathbf{v} \in X$ and $\mathbf{z} \in V^s$, there holds\vspace*{-0.75mm}
	\begin{align}\label{eq:fem:skew2} 
		{b}(\mathbf{v}, \mathbf{z}, \mathbf{z})
		= \tfrac{1}{2} (\{\mathrm{div}\,\Sigma_h \mathbf{v}\} \mathbf{z}, \mathbf{z})_{\Omega}\,.\\[-5.5mm]\notag
	\end{align}
\end{lemma}
\begin{proof}
	All the claims follow from Hölder's inequality, Lemma \ref{lem:fem:sigma_stability}, and integration-by-parts (which is possible because $\Sigma_h$ maps to $H(\operatorname{div};\Omega)$), similar to Lemma \ref{lem:fem:admissible}.
\end{proof}

\subsection{Discrete problem formulation} \label{sub:fem:discrete_problem}

\hspace{5mm}Since the boundary conditions fix volume flux everywhere on the topological boundary $\partial\Omega$, the combination of discretized boundary and divergence data have to fulfill a compatibility condition. Therefore, we define 
\begin{align}\label{def:g_1_tilde}
	g_1^h \coloneqq  g_1 + \langle\mathrm{div}\,\mathbf{g}_b^h - g_1\rangle_{\Omega} \in \smash{Y^{s'}} \,.
\end{align} 

We consider the following discrete counterpart of Problem (\hyperlink{Q}{Q}):

\textsc{Problem} (\hypertarget{Qh}{Q$_h$}). 
Find $(\mathbf{v}_h, q_h) \in X_h \times Q_h$ with $\mathbf{v}_h=\mathbf{g}_b^h$ in $\tr X_h$ such that for every $(\mathbf{z}_h,y_h) \in V_h\times Y_h$, there holds 
\begin{align}
	(\mathbf{S}(\mathbf{Dv}_h), \mathbf{Dz}_h)_{\Omega} + b(\mathbf{v}_h, \mathbf{v}_h,\mathbf{z}_h ) - (q_h, \mathrm{div}\,\mathbf{z}_h )_{\Omega} &= (\mathbf{f}, \mathbf{z}_h )_{\Omega}\,, \label{eq:fem:qh_main}  \\
	(\mathrm{div}\,\mathbf{v}_h,y_h)_{\Omega} &= (g_1^h,  y_h)_{\Omega}\,. \label{eq:fem:qh_div}  
\end{align}
\begin{remark}
	The divergence theorem implies that Problem (\hyperlink{Qh}{Q$_h$}) depends only on the values of $\mathbf{g}_b^h$ on Dirichlet boundary nodes.
	Furthermore, if $\mathbf{g}_b^h = \Pi_h^X\mathbf{g}_b$, then the properties of the Fortin interpolation operator yield that the boundary condition $\mathbf{v}_h=\mathbf{g}_b^h$ in $\tr X_h$  does not depend on the choice of the trace lift\linebreak $\mathbf{g}_b$. Discrete divergence preservation by $\Pi_h^X$ implies compatibility of $g_1$ and $\mathbf{g}_b^h$, so that~$g_1^h = g_1$~in~that~case.
\end{remark}

In order to establish the existence of discrete solutions to Problem (\hyperlink{Qh}{Q$_h$}), a discrete counterpart of Problem (\hyperlink{P}{P}) is useful. As a preparation, we construct a discrete extension of the boundary data which has the correct divergence.

\begin{lemma} \label{lem:fem:gh}
	There exists a discrete vector field $\mathbf{g}_h  \in X_h$ that satisfies 
	\begin{alignat}{2}
		(\mathrm{div}\,\mathbf{g}_h , y_h)_{\Omega} &= (g_1^h, y_h)_{\Omega} \quad &&\textrm{ for all } y_h \in Y_h\,, \label{eq:fem:gh_1} \\[-0.5mm]
		\mathbf{g}_h  &= \mathbf{g}_b^h &&\textrm{ in } \tr X_h\,. \label{eq:fem:gh_2}\\[-6.5mm]\notag
	\end{alignat}
	and $	\mathbf{g}_h  \to \mathbf{g} \coloneqq  \mathbf{g}_b + {\mathcal{B}} (g_1 - \mathrm{div}\,\mathbf{g}_b - \langle g_1 - \mathrm{div}\,\mathbf{g}_b\rangle_{\Omega})$ in $X^s $ $(h\to 0) $,
	where $\mathcal{B}\colon \smash{Q^{s'}} \to V^s$ denotes the Bogovski{\u{\i}} operator and the limit $\mathbf{g}\in X^s$ satisfies \eqref{eq:fem:div_eq}. 
\end{lemma}
\begin{proof}
	Since $g_1^h - \mathrm{div}\,\mathbf{g}_b^h \in Q^{s'}$, we choose $\mathbf{g}_h  \coloneqq  \mathbf{g}_b^h + \Pi_h^X {\mathcal{B}} (g_1^h - \mathrm{div}\,\mathbf{g}_b^h) \in X_h $.
	By Assumption \ref{assum:fem:projection_x2_weak}, we have that $\mathbf{g}_b^h \to \mathbf{g}_b$ in $X^s$ $(h\to 0)$.
	In consequence, the convergence and stability follow from continuity and linearity of $\Pi_h^X$ and the Bogovski{\u{\i}} operator.
\end{proof}

In order to establish the existence of a discrete velocity vector field in Problem (\hyperlink{Qh}{Q$_h$}), we prefer a problem formulation \textit{``hiding''} the discrete pressure. Taking the discrete lift $\mathbf{g}_h\in X_h$ from Lemma \ref{lem:fem:gh} and making the ansatz\vspace*{-1.5mm}
\begin{align}\label{eq:discrete_ansatz}
	\mathbf{u}_h \coloneqq  \mathbf{v}_h - \mathbf{g}_h \in V_{h,0}\,,\\[-6.5mm]\notag
\end{align}
leads to the following discrete counterpart of Problem (\hyperlink{P}{P}):

\textsc{Problem} (\hypertarget{Ph}{P$_h$}). 
Given a solution $\mathbf{g}_h   \in  X_h$ of \eqref{eq:fem:gh_1}, \eqref{eq:fem:gh_2}, find $\mathbf{u}_h \in  V_{h,0}$ such that for every $\mathbf{z}_h \in V_{h,0}$, there holds\vspace*{-1.5mm}
\begin{align*}
	(\mathbf{S}(\mathbf{Du}_h+\mathbf{Dg}_h), \mathbf{Dz}_h)_\Omega + b(\mathbf{u}_h+\mathbf{g}_h, \mathbf{u}_h+\mathbf{g}_h,\mathbf{z}_h)
	= (\mathbf{f}, \mathbf{z}_h)_\Omega\,.\\[-6.5mm]\notag
\end{align*} 

By Lemma \ref{lem:fem:gh} and  Lemma \ref{lem:fem:discrete-inf-sup}, Problem (\hyperlink{Qh}{Q$_h$}) and Problem (\hyperlink{Ph}{P$_h$}) are equivalent.

\subsection{Well-posedness and convergence analysis} \label{sub:fem:ex_conv}

\hspace*{5mm}Analogously to \cite[Prop.\ 3.3, Thm.\ 3.4]{JK23_fem_arxiv}, we establish the well-posedness (\textit{i.e.}, existence of discrete solutions), stability (\textit{i.e.}, \textit{a priori} bounds), and (weak) convergence of Problem (\hyperlink{Qh}{Q\textsubscript{$h$}}) (and Problem (\hyperlink{Ph}{P\textsubscript{$h$}}), respectively). While in \cite{JK23_fem_arxiv}, only the case $p \geq \tfrac{2d}{d+1}$ is treated, one easily checks that the argumentation works out indeed completely analogously, as it is the case in the continuous existence theory. Due to the length of the proofs, we refrain from repeating~them~here.

\begin{proposition} \label{prop:fem:existence}
	Let Assumptions \ref{assum:extra_stress}, \ref{assum:fem:data}, \ref{assum:fem:geometry}, \ref{assum:fem:projection_x}, \ref{assum:fem:projection_x2_weak} and \ref{assum:fem:reconstruction} be \mbox{fulfilled}.
	Moreover, we assume that $p>2$ or that the data are sufficiently small, i.e., $\norm{g_1}_s + \norm{\mathbf{g}_b}_{1,s} \leq C(\mathbf{S}, \mathbf{f}, \Omega)$.
	Then, there exists a solution $(\mathbf{v}_h, q_h)\in V_h\times Q_h$ to \mbox{Problem} (\hyperlink{Qh}{Q\textsubscript{h}}) and a constant $R>0$,  which depends on $\mathbf{S}$, $\delta$, $\norm{g_1}_{s,\Omega}$, $\norm{\mathbf{g}_b^h}_{1,s,\Omega}$, $\norm{\mathbf{f}}_{V^*}$, such that\vspace*{-3mm}
	\begin{align} \label{eq:fem:apriori}
		\|\mathbf{v}_h\|_{1,p,\Omega} + \|q_h\|_{s',\Omega} \leq R\,.
	\end{align} 
\end{proposition} 

\begin{proposition} \label{prop:fem:convergence}
	Let the assumptions of Proposition \ref{prop:fem:existence} be satisfied. If ${p \leq \tfrac{3d}{d+2}}$, then let,  in addition,  Assumption \ref{assum:fem:localbasis} (locally supported basis) be satisfied. Moreover, let $(h_n)_{n\in \mathbb{N}}\subseteq (0,1]$ be a sequence such that $h_n\to 0$ $(n\to \infty)$ and  
	let $(\mathbf{v}_{h_n}, q_{h_n})\in X_{h_n}\times Q_{h_n}$, $n\in \mathbb{N}$, be the corresponding sequence of solutions to Problem (\hyperlink{Qh}{Q$_{h_n}$}) which fulfill a uniform \textit{a priori}~estimate~\eqref{eq:fem:apriori}. 
	
	Then, there exists a subsequence $(n_k)_{k\in \mathbb{N}}\subseteq \mathbb{N}$ such that
	\begin{align*}
		\begin{aligned}
			\mathbf{v}_{h_{n_k}} &\rightharpoonup\mathbf{v} &&\quad\textrm{ in } X&&\quad(k\to \infty )\,, \\
			\smash{q_{h_{n_k}}} &\rightharpoonup q &&\quad\textrm{ in } Q^s&&\quad(k\to \infty )\,,
		\end{aligned}
	\end{align*}
	where $\smash{(\mathbf{v}, q)}\in X\times Q^s$ is a solution to Problem (\hyperlink{Q}{Q}) satisfying the \textit{a priori} bound \eqref{eq:fem:apriori}. 
\end{proposition} 

\section{\textit{A priori} error estimates}\label{sec:error_estimates}

\hspace{5mm}In this section, we derive \textit{a priori} error estimates for the approximation of a \textit{``regular''} solution $(\mathbf{v}, q)\in X\times Q^s$ of Problem (\hyperlink{Q}{Q}) by a solution $(\mathbf{v}_h, q_h)\in X_h\times Q_h$ of Problem (\hyperlink{Qh}{Q$_h$}).

Therefore, we define the auxiliary coefficients 
\begin{align*}
r &
\coloneqq    \min \{2, p\}\,,\\
\ell &
\coloneqq \max \{2, p, s \}
= \left\{
\begin{aligned}
&s = \smash{(\tfrac{p^*}{2})'} &&\textrm{if } p < \tfrac{4d}{d+4}\,, \\ 
&2 &&\textrm{if } p \in [\tfrac{4d}{d+4}, 2]\,, \\ 
&p &&\textrm{else\,.}
\end{aligned}
\right.
\end{align*}

\begin{assumption}[regularity] \label{assum:fem:regularity}
	We assume that $(\mathbf{v}, q)\in X\times Q^s$ is a solution of Problem (\hyperlink{Q}{Q}) with  $\mathbf{F}(\mathbf{Dv}) \in (W^{1,2}(\Omega))^{d\times d}$, $\mathbf{v} \in (W^{2,(r^*)'}(\Omega))^d$, ${\mathbf{g}_b \in (W^{k,\max\{\ell,\beta\}}(\Omega))^d}$ (\textit{cf}.\ Assumption \ref{assum:fem:projection_x2}), and $q \in W^{1,p'}(\Omega)$.
\end{assumption}

\begin{remark} \label{rem:fem:regularity}
	\begin{itemize}[noitemsep,topsep=2pt,leftmargin=!,labelwidth=\widthof{(iii)}]
		\item[(i)] 
		The regularity $\mathbf{F}(\mathbf{Dv}) \in (W^{1,2}(\Omega))^{d\times d}$ is natural for $p$-Laplace-type problems~(\textit{cf}.~\cite{gia-mod-86,giu1}). It has been proved in the space periodic setting (\textit{cf}.\ \cite{hugo-petr-rose}) or for regular~\mbox{domains}~if~${d\hspace{-0.1em}=\hspace{-0.1em}2}$~(\textit{cf}.~\cite{KMS2}).
		
		\item[(ii)] If $d=2$ or if $d=3$ and $p \geq \tfrac{4}{3}$, then $\mathbf{F}(\mathbf{Dv}) \in (W^{1,2}(\Omega))^{d\times d}$ implies that $\mathbf{v} \in (W^{2,(r^*)'}(\Omega))^d$ (\textit{cf}.\ \cite[Lem.\ 4.5]{BDR10_strongsol}).
		Moreover, the assumptions on $\mathbf{g}_b$ come into play only if ${\mathbf{g}_b^h \neq  \Pi_h^X \mathbf{g}_b}$. In other words, Assumption \ref{assum:fem:regularity} is in accordance with previous results (\textit{cf}. \cite{JK23_fem_arxiv}).
		
		\item[(iii)] Due to $\delta>0$, from $\mathbf{F}(\mathbf{Dv}) \in (W^{1,2}(\Omega))^{d\times d}$, it follows that $\mathbf{v}\in  (W^{2,r}(\Omega))^d$. Thus,  by $\mathbf{v}\in  (W^{2,r}(\Omega))^d$ or by $\mathbf{v}\in   (W^{2,(r^*)'}(\Omega))^d$ and Sobolev embeddings (for this, one easily checks that $r^*>d$ if $p>\frac{d}{2}$ and $((r^*)')^*> d$ if  $p\leq\frac{d}{2}$), there holds $\mathbf{v}\in (W^{1, d+\varepsilon}(\Omega))^d\hookrightarrow (L^\infty(\Omega))^d$~for~some~${\varepsilon>0}$. 
		
		\item[(iv)] If $p>2$ due to $\delta>0$, from the regularity $\mathbf{F}(\mathbf{Dv}) \in (W^{1,2}(\Omega))^{d\times d}$, it follows that $q\in W^{1,p'}(\Omega)$ if $\mathbf{f}\in (L^{p'}(\Omega))^d$ and $(\delta+\vert \mathbf{Dv}\vert)^{2-p}\vert \nabla q\vert^2\in L^1(\Omega)$ if $\mathbf{f}\in (L^2(\Omega))^d$ (\textit{cf}.~\mbox{\cite[Lem.~2.6]{KR23_ldg2}}).
	\end{itemize}
\end{remark}


\begin{proposition}[velocity error estimate up to \textit{`small'} pressure error] \label{prop:fem:error}
	Let the Assumptions \ref{assum:extra_stress}, \ref{assum:fem:data}, \ref{assum:fem:geometry}, \ref{assum:fem:projection_x}, \ref{assum:fem:projection_x2}, \ref{assum:fem:reconstruction} and \ref{assum:fem:regularity} be satisfied. Then, there exists a constant $c_0>0$, depending only on the characteristics of $\mathbf{S}$, $\delta^{-1}$, $\gamma_0$, $m$, $k$, and $\Omega$, such that from 
	\begin{align} \label{eq:fem:error_smallness}
		\|\mathbf{v}\|_{1,(\smash{\frac{r^*}{2}})',\Omega}
		\leq c_0\,,
	\end{align}
	for every $\zeta>0$, it follows that
	\begin{align}\begin{aligned} \label{eq:fem:error_vel}
			\|\mathbf{F}(\mathbf{Dv}_h) - \mathbf{F}(\mathbf{Dv})\|_{2,\Omega}^2
			&\leq c_1\, h^2  + c_2\, \rho_{(\phi_{\vert\mathbf{Dv}\vert})\d,\Omega}(h\, \nabla q) + \zeta \,\|q-q_h\|_{\ell',\Omega}^2\,,
	\end{aligned}\end{align}
	where the constant
	$c_1 >0$ depends on $\zeta$, on the characteristics of $\mathbf{S}$, $\delta^{-1}$, $c_0^{-1}$, $\gamma_0$, $m$, $k$, $\Omega$, $\|\nabla q\|_{p',\Omega}$, $\|\mathbf{v}\|_{2,(r^*)',\Omega}$, $\|\nabla\mathbf{F}(\mathbf{Dv})\|_{2,\Omega}$,  $\|\nabla^k \mathbf{g}_b\|_{\max\{\beta,\ell\},\Omega}$, $\|\mathbf{g}_b\|_{1,s,\Omega}$, and $\|g_1\|_{s,\Omega}$,  and the constant 
	$c_2 >0$ depends on the characteristics of $\mathbf{S}$ and $\gamma_0$.\enlargethispage{3.5mm}
\end{proposition}

\begin{proof}
	First, abbreviating
	\begin{align*}
		\begin{aligned}
			\mathbf{u}&\coloneqq \mathbf{v}-\mathbf{g}_b\in V\,,&&\quad\mathbf{u}_h\coloneqq \mathbf{v}_h-\mathbf{g}_b^h\in V_h\,,\\
			\mathbf{e}_h &\coloneqq  \mathbf{v}_h - \mathbf{v} \in X\,,&&\,\quad\mathbf{r}_h \coloneqq  \mathbf{u}_h   - \mathbf{u} \in V\,,
		\end{aligned}
	\end{align*}
	we arrive at the decomposition
	\begin{align}\label{eq:fem:decomp_eh}
		\mathbf{e}_h
		= \Pi_h^X \mathbf{r}_h + \{\Pi_h^X \mathbf{v} - \mathbf{v}\} + \{\mathbf{g}_b^h  - \Pi_h^X \mathbf{g}_b\} \quad \text{ in }X\,.
	\end{align} 
	
	Due to the decomposition  \eqref{eq:fem:decomp_eh}, the $\varepsilon$-Young inequality \eqref{eq:fem:orliczyoung} with $\psi\!=\!\phi_{\vert\mathbf{Dv}\vert}$, \eqref{eq:growth_S}, the approximation properties of $\Pi_h^X$ (\textit{cf}.\ \cite[Thms.\ 3.4, 5.1]{BBDR12_fem}), and Lemma \ref{lem:fem:projection_x2_oappr}, we have that
	\begin{align}\begin{aligned} \label{eq:fem:error_1}
			c \,\|\mathbf{F}(\mathbf{Dv}_h) - \mathbf{F}(\mathbf{Dv})\|_{2,\Omega}^2 
			&\leq (\mathbf{S}(\mathbf{Dv}_h) - \mathbf{S}(\mathbf{Dv}), \mathbf{De}_h)_\Omega \\ 
			&\leq (\mathbf{S}(\mathbf{Dv}_h) - \mathbf{S}(\mathbf{Dv}), \mathbf{D}\Pi_h^X \mathbf{r}_h)_\Omega \\
			&\quad + c_\varepsilon\, \|\mathbf{F}(\mathbf{Dv})-\mathbf{F}(\mathbf{D}\Pi_h^X \mathbf{v})\|_{2,\Omega}^2
			\\&\quad + c_\varepsilon\, \rho_{\phi_{\vert\mathbf{Dv}\vert},\Omega} (\mathbf{Dg}_b^h -\mathbf{D} \Pi_h^X \mathbf{g}_b)  \\		
			&\quad + \varepsilon \,\|\mathbf{F}(\mathbf{Dv}_h) - \mathbf{F}(\mathbf{Dv})\|_{2,\Omega}^2\\
			&\leq (\mathbf{S}(\mathbf{Dv}_h) - \mathbf{S}(\mathbf{Dv}), \mathbf{D}\Pi_h^X \mathbf{r}_h)_\Omega
			\\&\quad+c_{\varepsilon} \,h^2\, \big\{ \|\nabla \mathbf{F}(\mathbf{Dv})\|_{2,\Omega}^2 + \|\nabla^k \mathbf{g}_b\|_{\max\{2,p,\beta\}, \Omega}^2\big\}
			\\&\quad+ \varepsilon\, \|\mathbf{F}(\mathbf{Dv}_h) - \mathbf{F}(\mathbf{Dv})\|_{2,\Omega}^2\,.
	\end{aligned}\end{align}
	Subtracting the momentum equations in Problem (\hyperlink{Q}{Q}) and Problem  (\hyperlink{Qh}{Q$_h$}), for every ${\mathbf{z}_h \in V_h}$, yields the error equation
	\begin{align}
		\begin{aligned} \label{eq:fem:error_align}
			(\mathbf{S}(\mathbf{Dv}_h) - \mathbf{S}(\mathbf{Dv}), \mathbf{Dz}_h)_\Omega 
			&= ((q_h - q)\mathbf{I},  \mathbf{Dz}_h)_\Omega- (\mathbf{v} \otimes \mathbf{v}, \mathbf{Dz}_h)_\Omega - b(\mathbf{v}_h, \mathbf{v}_h, \mathbf{z}_h)\,.
	\end{aligned}\end{align} 
	Since $\mathbf{z}_h=\Pi_h^X \mathbf{r}_h \in V_h$ is an admissible test function in \eqref{eq:fem:error_align}, we deduce that
	\begin{align}
		\begin{aligned} \label{eq:fem:error_2}
			(\mathbf{S}(\mathbf{Dv}_h) - \mathbf{S}(\mathbf{Dv}), \mathbf{D}\Pi_h^X \mathbf{r}_h)_\Omega &= ((q_h - q)\mathbf{I}, \nabla\Pi_h^X \mathbf{r}_h)_\Omega\\&\quad - \big\{ (\mathbf{v} \otimes \mathbf{v}, \nabla \Pi_h^X \mathbf{r}_h)_\Omega + b(\mathbf{v}_h, \mathbf{v}_h, \Pi_h^X \mathbf{r}_h)\big\} \\
			&\eqqcolon  I^1_h + I_h^2\,.
		\end{aligned}
	\end{align}
	
	Therefore, we need to estimate $I_h^1$ and $I_h^2$:
	
	\emph{ad $I_h^1$.}
	With $(y_h, \operatorname{div} [\mathbf{e}_h +\mathbf{v} - \Pi_h^X \mathbf{v}])_{\Omega} = (y_h, \operatorname{div} [\mathbf{g}_b^h - \Pi_h^X \mathbf{g}_b])_\Omega $ for any $y_h \in Q_h$,
	the $\varepsilon$-Young inequality \eqref{eq:fem:orliczyoung} with $\psi=\varphi_{\vert \mathbf{Dv}\vert}$ and $\psi=\vert\cdot\vert^{\ell'}$, \eqref{eq:growth_S}, the approximation properties of $\Pi_h^X$ and $\Pi_h^Y$ (\textit{cf}.\ \cite[Thm. 4.6]{DR07_interpolation}), and \eqref{eq:fem:ph-gb}, we find that 
	\begin{align}\label{eq:Ih1}
		\begin{aligned} 
		I_h^1
		&=\inf_{\eta_h\in Q_h}{\big\{((\eta_h - q)\mathbf{I},  \mathbf{De}_h )_\Omega  + ((\eta_h - q)\mathbf{I},  \mathbf{Dv} - \mathbf{D}\Pi_h^X \mathbf{v})_\Omega\big\}}\\
		&\quad + ((q_h - q)\mathbf{I},  \mathbf{D}\Pi_h^X \mathbf{g}_b - \mathbf{Dg}_b^h )_\Omega
		\, + \, (\eta_h - q_h,  \operatorname{div} [\Pi_h^X \mathbf{g}_b - \mathbf{g}_b^h] )_\Omega \\
		&\leq c_{\varepsilon}\,  \inf_{\eta_h\in Q_h}{\big\{\rho_{(\phi_{\vert\mathbf{Dv}\vert})\d, \Omega}(\eta_h - q) \big\}}\\&\quad+\varepsilon\,\big\{\|\mathbf{F}(\mathbf{Dv}_h)-\mathbf{F}(\mathbf{Dv})\|_{2,\Omega}^2+\|\mathbf{F}(\mathbf{Dv})-\mathbf{F}(\mathbf{D}\Pi_h^X \mathbf{v})\|_{2,\Omega}^2\big\}
		\\&\quad + c_{\zeta}\, \|\mathbf{Dg}_b-\mathbf{D}\Pi_h^X \mathbf{g}_b \|_{\ell,\Omega}^2  + \zeta\, \|q_h - q\|_{\ell',\Omega}^2
		\\&\leq c_{\varepsilon}\, \big\{\rho_{(\phi_{\vert\mathbf{Dv}\vert})\d, \Omega}(h\,\nabla q)+h^2\,\|\nabla\mathbf{F}(\mathbf{Dv})\|_{2,\Omega}^2\big\}+c_\zeta\,\|\nabla^k \mathbf{g}_b\|_{\max \{\ell, \beta\},\Omega}^2
		\\&\quad +\varepsilon\,\|\mathbf{F}(\mathbf{Dv}_h)-\mathbf{F}(\mathbf{Dv})\|_{2,\Omega}^2+\zeta\, \|q_h - q\|_{\ell',\Omega}^2 \,.
	\end{aligned}
	\end{align}

	\emph{ad $I_h^2$.} If Temam's modification (\textit{cf}.\ \eqref{eq:fem:defi_conv_temam}) is used,  we proceed as in \cite[Thm. 4.2]{JK23_fem_arxiv}. Else if the reconstruction version (\textit{cf}.\ \eqref{eq:fem:defi_conv_recon}) is used, the consistency \eqref{eq:fem:conv-recon-consistency} yields that
	\begin{align}\begin{aligned} \label{eq:fem:decomp_j2-1}
		I_h^2
			&= - (\mathbf{v}\otimes (\mathbf{v} - \Sigma_h \mathbf{v}), \nabla\Pi_h^X \mathbf{r}_h)_{\Omega}+\big\{{b}(\mathbf{v}, \mathbf{v}, \Pi_h^X \mathbf{r}_h)-b(\mathbf{v}_h, \mathbf{v}_h, \Pi_h^X \mathbf{r}_h)\big\}\\&
			\eqqcolon I_h^{21}+I_h^{22} \,.
	\end{aligned}\end{align}
	
	First, using  Hölder's inequality, Definition \ref{defi:fem:reconstruction}(ii), Sobolev stability of $\Pi_h^X$ (\textit{cf}.\ \cite[Cor.\ 4.8]{DR07_interpolation}), the $\varepsilon$-Young inequality \eqref{eq:fem:orliczyoung} with $\psi=\vert\cdot\vert^2$, a Sobolev embedding (as $((r^*)')^*\ge r'$), and Korn's inequality, we find that
	\begin{align}\label{eq:fem:decomp_j2-1.2}
		\begin{aligned}
		\vert I_h^{21}\vert 
		&\leq \|\mathbf{v}\|_{\infty,\Omega} \|\mathbf{v}- \Sigma_h \mathbf{v}\|_{r',\Omega} \|\nabla\Pi_h^X \mathbf{r}_h\|_{r,\Omega} \\
		&\leq c_{\varepsilon}\, h^2\, \|\mathbf{v}\|_{\infty,\Omega}^2\|\mathbf{v}\|_{2, (r^*)',\Omega}^2 + \varepsilon\, \|\mathbf{Dr}_h\|_{r,\Omega}^2 \,.
	\end{aligned}
	\end{align}
	
	Second, with the trilinearity of ${b}$ and the identity $\Pi_h^X \mathbf{r}_h = \Pi_h^X \mathbf{e}_h + \Pi_h^X \mathbf{g}_b - \mathbf{g}_b^h$ in $X_h$, we decompose
	\begin{align}\begin{aligned} \label{eq:fem:decomp_j2-2}
			I_h^{22} 
			&= {b}(\mathbf{v}, \mathbf{v} - \Pi_h^X \mathbf{v}, \Pi_h^X \mathbf{r}_h)
				- {b}(\mathbf{e}_h, \Pi_h^X \mathbf{v}, \Pi_h^X \mathbf{r}_h)
			\\&\quad	- {b}(\mathbf{v}_h, \Pi_h^X \mathbf{r}_h, \Pi_h^X \mathbf{r}_h) 
			+ {b}(\mathbf{v}_h, \Pi_h^X \mathbf{g}_b - \mathbf{g}_b^h, \Pi_h^X \mathbf{r}_h) \\
			&\eqqcolon  J_h^1 +J_h^2 + J_h^3 + J_h^4 \,.
	\end{aligned}\end{align}

	So, it is left to estimate $J_h^{i}$, $i=1,\ldots,4$:
	
	\emph{ad $J_h^1$.} 
	Using Hölder's inequality, the $\varepsilon$-Young inequality \eqref{eq:fem:orliczyoung} with ${\psi=\vert\cdot\vert^2}$, Lemma \ref{lem:fem:sigma_stability} in connection with Remark \ref{rem:fem:regularity}(iii), the $L^{r'}$-approximability of $\Pi_h^X$ (\textit{cf}. \cite[Thm.\ 4.6]{DR07_interpolation}), the Sobolev embedding $W^{2,(r^*)'}(\Omega) \hookrightarrow  W^{1,r'}(\Omega)$ (one easily checks that $((r^*)')^*\ge r'$), stability of $\Pi_h^X$ and Korn's inequality, we obtain
	\begin{align}\label{eq:Ih21}
		\begin{aligned}
		\vert J^{1}_h\vert 
		&\leq \|\Sigma_h\mathbf{v}\|_{\infty,\Omega} \|\mathbf{v} - \Pi_h^X \mathbf{v}\|_{r',\Omega} \|\nabla\Pi_h^X \mathbf{r}_h\|_{r,\Omega}\\
		&\leq c_{\varepsilon} \, c(\|\mathbf{v}\|_{2,r,\Omega}, \|\mathbf{v}\|_{2, (r^*)',\Omega}) \, h^2 \|\mathbf{v}\|_{2, (r^*)',\Omega}^2
		+ \varepsilon\, \|\mathbf{Dr}_h\|_{r,\Omega}^2 \,.
	\end{aligned}
	\end{align}
	
	\emph{ad $J_h^{2}$.}
	For the next steps, we use the auxiliary function 
	\begin{align*} 
		\mathbf{s}_h \coloneqq \tfrac{1}{d} \langle\operatorname{div}\mathbf{g}_b^h - \operatorname{div}\mathbf{g}_b\rangle_{\Omega}\, \mathbf{id}_{\mathbb{R}^d} \in X_h\,,
	\end{align*}
	which satisfies
	\begin{align} \label{eq:fem:sh_div}
		\operatorname{div}\mathbf{s}_h=\langle\operatorname{div}\mathbf{g}_b^h - \operatorname{div}\mathbf{g}_b\rangle_{\Omega}\quad\text{ in }\Omega\,.
	\end{align} 
	For every $m \in \mathbb{N}$, $t \in [1, \infty]$, on the basis of Hölder's inequality and \eqref{eq:fem:proj_x2_sappr},  it holds  that
	\begin{align} \label{eq:fem:sh_estim}
		\begin{aligned} 
			\|\mathbf{s}_h\|_{m,t}
			&\leq c(\Omega,m,t)\, \|\nabla \mathbf{g}_b^h-\nabla\mathbf{g}_b\|_{1,\Omega} 
			\\&	\leq c(\Omega,m,t)\, h\, \|\nabla^k \mathbf{g}_b\|_{\beta,\Omega}\,.
		\end{aligned}
	\end{align}
	Moreover, from \eqref{eq:fem:q_div}, \eqref{eq:fem:qh_div}, \eqref{def:g_1_tilde}, and \eqref{eq:fem:sh_div}, it follows that $(\operatorname{div} (\mathbf{e}_h - \mathbf{s}_h), y_h)= 0$ for all $y_h \in Y_h$, so that, due to Definition \ref{defi:fem:reconstruction}(i), there holds
	\begin{align} \label{eq:fem:sh_Xh}
		\operatorname{div} \Sigma_h \mathbf{e}_h = \operatorname{div} \Sigma_h \mathbf{s}_h\,.
	\end{align}
	Using integration by parts (which is possible since $\Sigma_h \mathbf{v} \in H(\operatorname{div};\Omega)$), \eqref{eq:fem:sh_Xh}, Hölder's inequality, Definition \ref{defi:fem:reconstruction}(iii), \eqref{eq:fem:sh_estim}, Sobolev embeddings, the stability of $\Pi_h^X$ (\textit{cf}.\ \cite[Cor.\ 4.8]{DR07_interpolation}), Lemma \ref{lem:fem:sigma_stability}, the $\varepsilon$-Young inequality \eqref{eq:fem:orliczyoung} with ${\psi=\vert\cdot\vert^2}$, and Korn's inequality, we obtain
	\begin{align}\label{eq:Ih22}
		\begin{aligned} 
			\vert J^{2}_h\vert
			&= \vert (\operatorname{div} \Sigma_h \mathbf{e}_h, \Pi_h^X \mathbf{v} \cdot \Pi_h^X \mathbf{r}_h)_\Omega
			+ (\Pi_h^X \mathbf{r}_h \otimes \Sigma_h \mathbf{e}_h, \nabla\Pi_h^X \mathbf{v})_\Omega\vert  \\
			&\leq \|\operatorname{div} \Sigma_h \mathbf{s}_h\|_{r^*,\Omega} \|\Pi_h^X \mathbf{v}\|_{\smash{(\frac{r^*}{2})'},\Omega} \|\Pi_h^X \mathbf{r}_h\|_{r^*,\Omega} 
			\\&\quad+ \|\Pi_h^X \mathbf{r}_h\|_{r^*,\Omega} \|\Sigma_h \mathbf{e}_h\|_{r^*,\Omega} \|\nabla\Pi_h^X \mathbf{v}\|_{\smash{\smash{(\frac{r^*}{2})'}},\Omega} \\
			&\leq c_{\varepsilon} \, h^2\, \|\mathbf{v}\|_{\infty,\Omega}^2 \|\nabla^k \mathbf{g}_b\|_{\beta,\Omega}^2 + \varepsilon \,\|\mathbf{Dr}_h\|_{r,\Omega}^2
			\\&\quad+ c \,\|\mathbf{v}\|_{1,\smash{\smash{(\frac{r^*}{2})'}},\Omega} \|\mathbf{Dr}_h\|_{r,\Omega} \|\mathbf{e}_h\|_{1,r,\Omega}\,.
		\end{aligned}
	\end{align}
	
	\emph{ad $J_h^{3}$.}
	According to \eqref{eq:fem:skew2}, Hölder's inequality, $\mathbf{v}_h = \mathbf{e}_h-\mathbf{s}_h + \mathbf{s}_h + \mathbf{v}$, \eqref{eq:fem:sh_Xh}, Definition \ref{defi:fem:reconstruction}(i),(iii), a Sobolev embedding, Hölder's inequality and \eqref{eq:fem:sh_estim}, the $\varepsilon$-Young inequality \eqref{eq:fem:orliczyoung} with ${\psi=\vert\cdot\vert^2}$, and \textit{a priori} boundedness of $\|\mathbf{r}_h\|_{1,r,\Omega}$ (\textit{cf}.\ Proposition \ref{prop:fem:convergence}), there holds
	\begin{align}\label{eq:Ih23}
		\begin{aligned} 
		\vert J_h^{3}\vert
		&\leq  \|\operatorname{div} \Sigma_h \mathbf{v}_h\|_{\smash{\smash{(\frac{r^*}{2})'}},\Omega} \|\Pi_h^X\mathbf{r}_h\|_{r^*,\Omega}^2 \\
		&\leq c\, \|\operatorname{div} \Sigma_h (\mathbf{s}_h + \mathbf{v})\|_{\smash{\smash{(\frac{r^*}{2})'}},\Omega} \|\Pi_h^X\mathbf{r}_h\|_{1,r,\Omega}^2 \\
		&\leq c\, \big\{ h\, \|\nabla^k\mathbf{g}_b\|_{\beta,\Omega} + \|\mathbf{v}\|_{1,\smash{\smash{(\frac{r^*}{2})'}},\Omega} \big\} \|\mathbf{Dr}_h\|_{r,\Omega}^2 \\
		&\leq c \, \big\{ \varepsilon + \|\mathbf{v}\|_{1,\smash{\smash{(\frac{r^*}{2})'}},\Omega} \big\} \|\mathbf{Dr}_h\|_{r,\Omega}^2 + c_{\varepsilon}\, h^2\, \|\nabla^k\mathbf{g}_b\|_{\beta,\Omega}^2 \,.
		\end{aligned}
	\end{align}

	\emph{ad $J_h^{4}$.}
	Using integration by parts, Hölder's inequality, $\mathbf{v}_h = \mathbf{e}_h-\mathbf{s}_h + \mathbf{s}_h + \mathbf{v}$, \eqref{eq:fem:sh_Xh}, Definition \ref{defi:fem:reconstruction}(i),(iii), a Sobolev embedding, Lemma \ref{lem:fem:sigma_stability}, stability of $\Pi_h^X$ (\textit{cf}.\ \cite[Cor.\ 4.8]{DR07_interpolation}), \eqref{eq:fem:ph-gb}, the $\varepsilon$-Young inequality \eqref{eq:fem:orliczyoung} with ${\psi=\vert\cdot\vert^2}$, and \eqref{eq:fem:apriori}, we find that
	\begin{align}\label{eq:Ih24}
		\begin{aligned}
		\vert J_h^{4}\vert
		&\leq \|\operatorname{div} \Sigma_h \mathbf{v}_h\|_{r^*,\Omega} \|\Pi_h^X \mathbf{g}_b - \mathbf{g}_b^h\|_{\smash{\smash{(\frac{r^*}{2})'}},\Omega} \|\mathbf{r}_h\|_{r^*,\Omega} \\
		&\quad + \|\Pi_h^X\mathbf{r}_h \|_{r^*,\Omega} \|\Sigma_h \mathbf{v}_h\|_{r^*,\Omega} \|\Pi_h^X \mathbf{g}_b - \mathbf{g}_b^h\|_{1, \smash{(\frac{r^*}{2})'},\Omega} \\
		&\leq c\, \|\operatorname{div} \Sigma_h (\mathbf{s}_h+\mathbf{v})\|_{r^*,\Omega} \|\Pi_h^X \mathbf{g}_b - \mathbf{g}_b^h\|_{1,\smash{(\frac{r^*}{2})'},\Omega} \|\mathbf{Dr}_h\|_{r,\Omega} \\
		&\quad + c\, \| \mathbf{r}_h \|_{1,r,\Omega} \|\mathbf{v}_h\|_{1,r,\Omega} \|\Pi_h^X \mathbf{g}_b - \mathbf{g}_b^h\|_{1, \smash{(\frac{r^*}{2})'},\Omega} \\
		&\leq c\,h^2\,\big\{1+c_\varepsilon\,\|\mathbf{v}\|_{2,r,\Omega}^2\big\} \,\|\nabla^k\mathbf{g}_b\|_{\max \{\ell,\beta \},\Omega}^2 
		+ \varepsilon\, \|\mathbf{Dr}_h\|_{r,\Omega}^2  \,.
	\end{aligned}
	\end{align}
	
	Due to $\mathbf{e}_h = \mathbf{r}_h + \mathbf{g}_b^h - \mathbf{g}_b$, Poincaré's inequality, Korn's inequality, \eqref{eq:fem:proj_x2_sappr}, \cite[Eq.\ (4.17)]{JK23_fem_arxiv} and the regularity of $\mathbf{v}$ imply
	\begin{align}\label{eq:fem:error_12} 
		\begin{aligned} 
		\|\mathbf{e}_h\|_{1,r,\Omega}^2+ \|\mathbf{r}_h\|_{1,r,\Omega} ^2  
		&\leq c\, \|\mathbf{De}_h\|_{r,\Omega}^2 + c\, h^2\, \|\nabla^k \mathbf{g}_b\|_{\max \{r, \beta\},\Omega}^2\\
		&	\leq c\,  
		\|\mathbf{F}(\mathbf{Dv}_h) - \mathbf{F}(\mathbf{Dv})\|_{2,\Omega}^2
		+ c\, h^2\, \|\nabla^k \mathbf{g}_b\|_{\max \{r, \beta\},\Omega}^2\,. 
	\end{aligned}
	\end{align} 
	 
	As a consequence, if we combine \eqref{eq:fem:decomp_j2-1}--\eqref{eq:fem:error_12}, we find that
	\begin{align}\label{eq:fem:error_13}
		\vert I_h^2\vert \leq c_1\, h^2+ c_2 \,\big\{\varepsilon + \|\mathbf{v}\|_{1,\smash{(\frac{r^*}{2})'},\Omega}\big\}\, \|\mathbf{F}(\mathbf{Dv}_h) - \mathbf{F}(\mathbf{Dv})\|_{2,\Omega}^2\,.
	\end{align}
	
	Eventually, if we use \eqref{eq:Ih1} and \eqref{eq:fem:error_13} in \eqref{eq:fem:error_2}, from \eqref{eq:fem:error_1}, we deduce that
	\begin{align*}
		\|\mathbf{F}(\mathbf{Dv}_h) - \mathbf{F}(\mathbf{Dv})\|_{2,\Omega}^2
		&\leq c\,\big\{\varepsilon + \|\mathbf{v}\|_{1,\smash{(\frac{r^*}{2})'},\Omega} \big\}\, \|\mathbf{F}(\mathbf{Dv}_h) - \mathbf{F}(\mathbf{Dv})\|_{2,\Omega}^2 \\
		&\quad + c\,(c_{\zeta} + c_\varepsilon)\, h^2 
		+ c_{\varepsilon}\, \rho_{(\phi_{\vert\mathbf{Dv}\vert})\d,\Omega}(h\, \nabla q) + \zeta\, \|q-q_h\|_{\ell',\Omega}^2 \,.
	\end{align*}
	Choosing $\varepsilon$ sufficiently small and employing the smallness assumption \eqref{eq:fem:error_smallness} allows to absorb the error term from the right hand side, which yields the claim \eqref{eq:fem:error_vel}.
\end{proof} 

The next step is to bound the pressure error up to the velocity error.

\begin{proposition}[pressure error estimate up to the velocity error] \label{prop:fem:error_pres}
	Suppose the Assumptions \ref{assum:extra_stress}, \ref{assum:fem:data}, \ref{assum:fem:geometry}, \ref{assum:fem:projection_x}, \ref{assum:fem:projection_x2}, \ref{assum:fem:reconstruction} and \ref{assum:fem:regularity}.
	Then, there holds
	\begin{align} \label{eq:fem:error_pres}
		&	\begin{aligned} 
			\norm{q_h - q}_{s',\Omega}
			&\leq c_1 \,h \,\big\{\|\nabla \mathbf{F}(\mathbf{Dv})\|_{2,\Omega}+\|\nabla q\|_{s',\Omega}\big\}\\&\quad + c_2\, \|\mathbf{F}(\mathbf{Dv}_h) - \mathbf{F}(\mathbf{Dv})\|_2^{\smash{\min \{1, \tfrac{2}{p'}\}}} \,,
		\end{aligned} \\
		&\begin{aligned} 
			\|q_h - q\|_{\ell',\Omega}
			&\leq c_1 \,h\,\big\{\|\nabla \mathbf{F}(\mathbf{Dv})\|_{2,\Omega}+\|\nabla q\|_{\ell',\Omega}\big\}\\&\quad  + c_2 \,\|\mathbf{F}(\mathbf{Dv}_h) - \mathbf{F}(\mathbf{Dv})\|_{2,\Omega}\,,
		\end{aligned}  \label{eq:fem:2error_pres}
	\end{align}
	where the constant $c_1>0$ depends only on $\|\mathbf{v}\|_{1,p,\Omega}$ and $\gamma_0$, and the constant $c_2>0$ depends only on the characteristics of $\mathbf{S}$, $\delta$, and $ \|\mathbf{v}\|_{1,p,\Omega}$. 
\end{proposition}

The proof of Proposition \ref{prop:fem:error_pres} follows the same reasoning as in \cite[Thm. 4.4, Cor. 4.7]{JK23_fem_arxiv} using the following error estimate for the convective term.

\begin{lemma} \label{lem:fem:estim_conv_pressure}
	Suppose the assumptions from Proposition \ref{prop:fem:error_pres}.
	Then, for every $\mathbf{v}_h \in X_h$,~there~holds
	\begin{align*} 
			\vert b(\mathbf{v}_h, \mathbf{v}_h, \mathbf{z}_h)
				+ (\mathbf{v} \otimes \mathbf{v}, \nabla\mathbf{z}_h)_{\Omega}\vert
			\leq c\,\big\{h\, \|\mathbf{v}\|_{2,r,\Omega} + \|\mathbf{e}_h\|_{1,r,\Omega}\big\}\, \|\mathbf{z}_h\|_{1,s,\Omega}\,,
	\end{align*}
	where the constant $c>0$ depends on $\gamma_0$, $\delta$, $\mathbf{S}$, $ \| \mathbf{F}(\mathbf{Dv}) \|_{1,2,\Omega}$, $ \| \mathbf{g}_b \|_{1,s,\Omega}$ and $\norm{q}_{1,p',\Omega}$. 
\end{lemma}
\begin{proof}
	If Temam's modification \eqref{eq:fem:defi_conv_temam} is used as discrete convective term, the estimate is proved analogously to \cite[Thm. 4.4]{JK23_fem_arxiv}. 
		For the convenience of the reader, we briefly recapitulate the key arguments here.
		Using consistency \eqref{eq:fem:temam_consistent}, we obtain the decomposition
		\begin{align*}
			&b(\mathbf{v}_h, \mathbf{v}_h, \mathbf{z}_h)
			+ (\mathbf{v} \otimes \mathbf{v}, \nabla\mathbf{z}_h)_{\Omega} \\
			&= \tfrac{1}{2} (g_1 (\mathbf{v}_h-\mathbf{v}), \mathbf{z}_h) _\Omega
			- \tilde{b}(\mathbf{v}, \mathbf{v} - \Pi_h^X \mathbf{v}, \mathbf{z}_h)
			+ \tilde{b}(\mathbf{e}_h, \Pi_h^X \mathbf{v}, \mathbf{z}_h)
			+ \tilde{b}(\mathbf{v}_h, \Pi_h^X\!\mathbf{e}_h, \mathbf{z}_h)
			\\&	\eqqcolon I_h^2 + I_h^{31}+I_h^{32}+I_h^{33} \,.
		\end{align*}
		With Hölder's inequality, Sobolev embeddings and the stability and approximation properties of $\Pi_h^X$, we now estimate
		\begin{align*}
			\vert I_h^2\vert 
			&\leq c \, \|g_1\|_{s,\Omega} \|\mathbf{e}_h\|_{r^*,\Omega}\|\mathbf{z}_h\|_{p^*,\Omega} \\
			&\leq c \, \|g_1\|_{s,\Omega} \|\mathbf{e}_h\|_{1,r,\Omega}\|\mathbf{z}_h\|_{1,p,\Omega} \,,\\
			\vert I_h^{31}\vert 
			&\leq \|\mathbf{v}\|_{\infty,\Omega}  \|\mathbf{z}_h\|_{1,(r^*)',\Omega} \|\mathbf{v} - \Pi_h^X \mathbf{v}\|_{1,r,\Omega} 
			\\&\leq c \,h\,\|\mathbf{v}\|_{\infty,\Omega}  \|\mathbf{z}_h\|_{1,s,\Omega} \|\nabla^2 \mathbf{v}\|_{r,\Omega} \,, \\
			\vert I_h^{32}\vert 
			&\leq c\, \|\mathbf{e}_h\|_{r^*,\Omega} \|\mathbf{z}_h\|_{1,s,\Omega} \|\Pi_h^X \mathbf{v}\|_{1, r^*,\Omega} \\
			&\leq c\, \|\mathbf{e}_h\|_{1,r,\Omega} \|\mathbf{z}_h\|_{1,s,\Omega} \| \mathbf{v}\|_{1, r^*,\Omega} \,, \\
			\vert I_h^{33}\vert 
			&\leq c\, \|\mathbf{v}_h\|_{p^*,\Omega} \|\mathbf{z}_h\|_{1,s,\Omega} \|\Pi_h^X \mathbf{e}_h\|_{1,p,\Omega} \\
			&\leq c\, \|\mathbf{v}_h\|_{1,p,\Omega} \|\mathbf{z}_h\|_{1,s,\Omega} \| \mathbf{e}_h\|_{1,p,\Omega} \,.
		\end{align*}
		Adding these inequalities yields the result.
	
	Else, if the reconstruction version \eqref{eq:fem:defi_conv_recon} is employed, using consistency \eqref{eq:fem:conv-recon-consistency}, we decompose as follows:
	\begin{align}\label{eq:fem:estim_conv_pressure.1}
		\begin{aligned} 
			b(\mathbf{v}_h, \mathbf{v}_h, \mathbf{z}_h)
			+ (\mathbf{v} \otimes \mathbf{v}, \nabla\mathbf{z}_h)_{\Omega} 
			&= - b(\mathbf{v}, \mathbf{v} - \Pi_h^X \mathbf{v}, \mathbf{z}_h)
			+ b(\mathbf{e}_h, \Pi_h^X \mathbf{v}, \mathbf{z}_h) \\&\quad + b(\mathbf{v}_h, \Pi_h^X \mathbf{e}_h, \mathbf{z}_h)
			  + (\mathbf{v} \otimes \{\mathbf{v} - \Sigma_h \mathbf{v}\}, \nabla\mathbf{z}_h)_{\Omega} \\
			&\eqqcolon  I_h^{1} + I_h^{2} + I_h^{3} + I_h^{4}\,.
		\end{aligned}
	\end{align}
	So, it is left to estimate $I_h^{i} $, $i=1,\ldots,4$:
	
	\emph{ad $I_h^{1}$.} Using  Hölder's inequality, Lemma \ref{lem:fem:sigma_stability}, a Sobolev embedding, and the approximation properties of $\Pi_h^X$ (\textit{cf}.\ \cite[Thm.\ 4.6]{DR07_interpolation}), we find that
	\begin{align}\label{eq:fem:estim_conv_pressure.2}
		\begin{aligned}
		\vert I_h^{1}\vert 
		&\leq \|\Sigma_h \mathbf{v}\|_{r^*,\Omega} \|\mathbf{v} - \Pi_h^X \mathbf{v}\|_{r^*,\Omega} \|\nabla\mathbf{z}_h\|_{\smash{\smash{(\frac{r^*}{2})'}},\Omega}
		\\&\leq c\,h\, \|\mathbf{v}\|_{1,r,\Omega} \|\mathbf{v}\|_{2,r,\Omega} \|\mathbf{z}_h\|_{1,s,\Omega} \,.
		\end{aligned}
	\end{align}
	
	\emph{ad $I_h^{2}$ \& $I_h^{3}$.}
	Using  Hölder's inequality, Lemma \ref{lem:fem:sigma_stability}, a Sobolev embedding, and the stability properties of $\Pi_h^X$ \cite[Thm. 4.5]{DR07_interpolation}, we obtain
	\begin{align}\label{eq:fem:estim_conv_pressure.3}
		\begin{aligned}
		\vert I_h^{2}\vert 
		&\leq c\, \|\mathbf{e}_h\|_{1,r,\Omega} \|\mathbf{v}\|_{1,r,\Omega} \|\mathbf{z}_h\|_{1,s,\Omega}\,, \\
		\vert I_h^{3}\vert 
		&\leq c\, \|\mathbf{e}_h\|_{1,r,\Omega} \|\mathbf{v}_h\|_{1,r,\Omega} \|\mathbf{z}_h\|_{1,s,\Omega} \,,
		\end{aligned}
	\end{align}
	where $\|\mathbf{v}_h\|_{1,r,\Omega}$ is bounded by given data (\textit{cf}. \eqref{eq:fem:apriori}) and, hence, by $\delta$, $\mathbf{S}$ and the regularity of $\mathbf{v}$, $\mathbf{g}_b$, $q$ (\textit{cf.}\ Assumption \ref{assum:fem:regularity}).
	
	\emph{ad $I_h^{4}$.} Using Hölder's inequality, a Sobolev embedding, and consistency of $\Sigma_h$ (\textit{cf}.\ Definition \ref{defi:fem:reconstruction}(ii)), we deduce that
	\begin{align}\label{eq:fem:estim_conv_pressure.4}
		\vert I_h^{4}\vert
		\leq c\,h\, \|\mathbf{v}\|_{1,r,\Omega} \|\mathbf{v}\|_{2,r,\Omega} \|\mathbf{z}_h\|_{1,s,\Omega} \,.
	\end{align}
	Eventually, combining \eqref{eq:fem:estim_conv_pressure.2}--\eqref{eq:fem:estim_conv_pressure.4} in \eqref{eq:fem:estim_conv_pressure.1}, we conclude the assertion.
\end{proof}

\begin{proof}[Proof (of Proposition \ref{prop:fem:error_pres})]
	By means of  Lemma \ref{lem:fem:estim_conv_pressure}, the assertion follows by the same reasoning as in \cite[Thm.\ 4.4, Cor.\ 4.7]{JK23_fem_arxiv}.
\end{proof}

\begin{theorem}[velocity error estimate] \label{thm:fem:error_vel_final}
	Under the Assumptions \ref{assum:extra_stress}, \ref{assum:fem:data}, \ref{assum:fem:geometry}, \ref{assum:fem:projection_x}, \ref{assum:fem:projection_x2}, \ref{assum:fem:reconstruction} and \ref{assum:fem:regularity}, there holds 
	\begin{align} \label{cor:fem:error_vel.1}
		\begin{aligned}
			\|\mathbf{F}(\mathbf{Dv}_h) - \mathbf{F}(\mathbf{Dv})\|_{2,\Omega}^2 
			&\leq c_1\, h^2+ c_2\, \rho_{(\phi_{\vert\mathbf{Dv}\vert})\d,\Omega}(h\, \nabla q) \\
			&\leq c_1 \,h^2  + c_2\, h^{\smash{\min\{2,p'\}}}\,\rho_{(\phi_{\vert\mathbf{Dv}\vert})\d,\Omega}(\nabla q)\,,
	\end{aligned}\end{align}
	where the constant
	$c_1 >0$ depends on the characteristics of $\mathbf{S}$, $\delta^{-1}$, $\gamma_0$, $m$, $k$, $\Omega$, $\|\nabla q\|_{p',\Omega}$, $\|\mathbf{v}\|_{2,(r^*)',\Omega}$, $\|\nabla\mathbf{F}(\mathbf{Dv})\|_{2,\Omega}$,  $\|\nabla^k \mathbf{g}_b\|_{\max\{\beta,\ell\},\Omega}$, $\|\mathbf{g}_b\|_{1,s,\Omega}$, and $\|g_1\|_{s,\Omega}$,  and the constant 
	$c_2 >0$ depends on the characteristics of $\mathbf{S}$ and $\gamma_0$. 
	If $p > 2$ and, in addition, $\mathbf{f} \in (L^2(\Omega))^d$,~then~there~holds
	\begin{align}
		\label{cor:fem:error_vel.2}
		\begin{aligned}
			\|\mathbf{F}(\mathbf{Dv}_h) - \mathbf{F}(\mathbf{Dv})\|_{2,\Omega}^2 
			\leq c_1 \, h^2  +c_2 \,h^2\,\|(\delta+\vert \mathbf{Dv}\vert)^{\smash{\frac{p-2}{2}}}\nabla q\|_{2,\Omega}^2 \,.
		\end{aligned}
	\end{align}
\end{theorem}
\begin{proof}
	Inserting \eqref{eq:fem:2error_pres} into \eqref{eq:fem:error_vel} and 
for $\zeta > 0$ sufficiently small,  
	we infer \eqref{cor:fem:error_vel.1}$_1$. Then, 
	arguing as in \cite{BBDR12_fem, KR23_ldg3}, from \eqref{cor:fem:error_vel.1}$_1$,  we conclude  \eqref{cor:fem:error_vel.1}$_2$ and \eqref{cor:fem:error_vel.2}.
\end{proof}

If 
$\{\mathcal{T}_h\}_{h>0}$  is quasi-uniform, a stronger stability statement follows from inverse estimates and the regularity of the velocity vector field (\textit{cf}.\ Assumption \ref{assum:fem:regularity}).

\begin{lemma}[improved stability for the discrete velocity] \label{lem:fem:strong_stability}
	Let $p < 2$, $\delta > 0$, and assume that $\{\mathcal{T}_h\}_{h>0}$  is quasi-uniform. Then,
	from the assumptions of Theorem \ref{thm:fem:error_vel_final}, for every $\alpha < \infty$ if $d=2$ and $\alpha < 3p$ if $d=3$, it follows that
	\begin{equation*}
		\|\mathbf{v}_h\|_{1,\alpha,\Omega} + \|\mathbf{e}_h\|_{1,\alpha,\Omega} \leq c\,,
	\end{equation*}
	where $c>0$ depends only on the same quantities as $c_1$ in Theorem \ref{thm:fem:error_vel_final}.
\end{lemma}

\begin{proof} 
	First, by Proposition \ref{prop:fem:existence}, for $\alpha_1 \coloneqq p$, we have that $\|\mathbf{v}_h\|_{1,\alpha_1,\Omega} < c_1$. Now, for $n\in \mathbb{N}$, we assume that $\|\mathbf{Dv}_h\|_{\alpha_n,\Omega} < c_1$ for some $\alpha_n \in (1, \alpha_{lim})$, where $\alpha_{lim} \coloneqq \infty$ if $d=2$ and $\alpha_{lim} \coloneqq 3p$ if $d=3$. Then,
	\cite[Lem.\ 4.1]{BDR10_strongsol} and Theorem \ref{thm:fem:error_vel_final} yield that
	\begin{equation} \label{eq:fem:extra_stability1}
		\|\mathbf{De}_h\|_{\beta_n, \Omega}
		\leq c\, \|\mathbf{F}(\mathbf{Dv}_h) - \mathbf{F}(\mathbf{Dv})\|_{2,\Omega} \|(\delta + \vert \mathbf{Dv}\vert + \vert \mathbf{Dv}_h\vert)^{2-p}\|_{\smash{\frac{\beta_n}{2-\beta_n}},\Omega}^{\textcolor{purple}{\smash{\frac{1}{2}}}}
		\leq c\, h\,,
	\end{equation}
	where $\beta_n \coloneqq \smash{\tfrac{2\alpha_n}{2-p+\alpha_n}}  \in (1,2)$ and, hence, $\tfrac{(2-p) \beta_n}{2-\beta_n} = \alpha_n$. From the identities $\mathbf{v}_h = \Pi_h^X \mathbf{e}_h + \Pi_h^X \mathbf{v}$, $\mathbf{e}_h = \Pi_h^X \mathbf{e}_h + \Pi_h^X \mathbf{v} - \mathbf{v}$, a global inverse estimate (\textit{cf}.\ \cite[Lem.\ 3.5]{Bartels15}), stability of $\Pi_h^X$ (\emph{cf.} \cite[Thm.\ 4.5]{DR07_interpolation}), \cite[Lem.\ 4.5]{BDR10_strongsol}, the embedding $W^{2,\smash{\frac{3p}{p+1}}}(\Omega) \hookrightarrow  W^{1,3p}(\Omega)$ if $d=3$ (if $d=2$, a similar argument yields that ${\mathbf{v} \in (W^{1,t}(\Omega))^d}$ for all $t < \infty$), and \eqref{eq:fem:extra_stability1}, we deduce that 
	\begin{align*}
		\|\mathbf{Dv}_h\|_{\alpha_{n+1},\Omega} + \|\mathbf{De}_h\|_{\alpha_{n+1},\Omega}
		&\leq 2\, \|\mathbf{D}\Pi_h^X \mathbf{e}_h\|_{\alpha_{n+1},\Omega} + c \, \|\mathbf{D} \mathbf{v}\|_{\alpha_{n+1},\Omega} \\
		&\leq c\, h^{\smash{\frac{d}{\alpha_{n+1}} - \frac{d}{\beta_n}}} \|\mathbf{D}\Pi_h^X \mathbf{e}_h\|_{\beta_n,\Omega} + c
		\leq c_1\,,
	\end{align*}
	if $\smash{\frac{d}{\alpha_{n+1}} - \frac{d}{\beta_n}} + 1 =0$ and $\alpha_{n+1} < \alpha_{lim}$. By definition, for every $n\in\mathbb{N}$, we have that
	\begin{equation} \label{eq:fem:extra_stability2}
		\alpha_{n+1} \coloneqq \frac{d\beta_n}{d-\beta_n}= \frac{2d \alpha_n}{2d - dp + \alpha_nd - 2 \alpha_n}=\begin{cases}
			\frac{2}{2-p}\alpha_n&\text{ if }d=2\,,\\
		 \tfrac{6 \alpha_n}{6-3p+\alpha_n}&\text{ if }d=3\,.
		\end{cases}
	\end{equation} 
	Next, we distinguish the cases $d=2$ and $d=3$:
	
	\emph{ad $d=2$.} From \eqref{eq:fem:extra_stability2}, due to $p<2$ and $\alpha_1=p$, it follows that $\alpha_n\to \alpha_{lim}$ $(n\to \infty)$.
	
	\emph{ad $d=3$.}
	If we set $\delta_n\coloneqq \alpha_{n} - 3p$ for $n\in \mathbb{N}$,  due to \eqref{eq:fem:extra_stability2} and $6+\delta_n >6-3p > 0$~for~all~$n\in \mathbb{N}$, we have that\enlargethispage{7.5mm} 
	\begin{equation} \label{eq:fem:extra_stability_3}
		\frac{\delta_{n+1}}{\delta_n} = \frac{6-3p}{6+\delta_n}\in (0,1)\,.
	\end{equation}
	Since $\delta_0 
	 < 0$, \eqref{eq:fem:extra_stability_3} implies that $\delta_n \in (\delta_0, 0)$ for all $n \in \mathbb{N}$. As a result, 
	 ${\tfrac{\delta_{n+1}}{\delta_n}
	  < \tfrac{6-3p}{6+\delta_0} < 1}$~for~all~$n \in \mathbb{N}$, which implies that $\delta_n\to 0$ $(n\to \infty)$ and, thus, $\alpha_n \to \alpha_{lim}$ $(n \to \infty)$.
	 
	 In summary,  for every $\alpha < \alpha_{lim}$, we find that
	\begin{equation} \label{eq:fem:extra_stability4}
		\|\mathbf{Dv}_h\|_{\alpha,\Omega} + \|\mathbf{De}_h\|_{\alpha,\Omega} < c_1\,.
	\end{equation}
 From a global inverse estimate (\textit{cf}.\ \cite[Lem.\ 3.5]{Bartels15}), Sobolev stability and approximability of $\Pi_h^X$ (\emph{cf.} \cite[Thms. 4.5, 4.6]{DR07_interpolation}) and Assumption \ref{assum:fem:projection_x2}(iii), we deduce that\vspace*{-0.5mm}
	\begin{align}\label{eq:fem:extra_stability5}
		\begin{aligned} 
			\|\mathbf{g}_b^h - \mathbf{g}_b\|_{1,\ell^*,\Omega}
			&\leq c \,h^{-1} \, \|\Pi_h^X (\mathbf{g}_b^h - \mathbf{g}_b)\|_{1,\ell,\Omega} + c\, \|\Pi_h^X \mathbf{g}_b - \mathbf{g}_b\|_{1,\ell^*,\Omega} \\[-0.5mm]
			&\leq c\, \|\mathbf{g}_b\|_{k, \max\{\ell, \beta\},\Omega}\,.
		\end{aligned}
	\end{align}
	Similarly to the first step in \eqref{eq:fem:error_12}, \eqref{eq:fem:extra_stability5}, Assumption \ref{assum:fem:projection_x2}, $\ell^* \geq \alpha_{lim}$, and \eqref{eq:fem:extra_stability4} imply\vspace*{-0.5mm}
	\begin{align*}
		\|\mathbf{v}_h\|_{1,\alpha,\Omega}+	\|\mathbf{e}_h\|_{1,\alpha,\Omega}
		&\leq 2	\|\mathbf{e}_h\|_{1,\alpha,\Omega}+	\|\mathbf{v}\|_{1,\alpha,\Omega}
		\\&\leq c\, (\|\mathbf{De}_h\|_{\alpha,\Omega} + \|\mathbf{g}_b^h - \mathbf{g}_b\|_{1,\alpha,\Omega}+\|\mathbf{v}\|_{1,\alpha,\Omega})
		\leq c\,,
	\end{align*}
	where we used that $\mathbf{v} \in (W^{1,3p}(\Omega))^d$ if $d=3$ and  $\mathbf{v} \in (W^{1,t}(\Omega))$ for $t < \infty$ if $d=2$.
\end{proof}

\begin{proposition}[Pressure error estimates]\label{prop:fem:error_pres2}
	Let the Assumptions \ref{assum:extra_stress}, \ref{assum:fem:data}, \ref{assum:fem:geometry}, \ref{assum:fem:projection_x}, \ref{assum:fem:projection_x2}, \ref{assum:fem:reconstruction} and \ref{assum:fem:regularity} be satisfied. Then, if
	$p < 2$ and $\{\mathcal{T}_h\}_{h>0}$ is quasi-uniform, there holds\vspace*{-0.5mm}
	\begin{align} 
			\|q_h - q\|_{2,\Omega}
			&\leq c \,h\,, \label{eq:fem:3error_pres} \\
			\|q_h - q\|_{p',\Omega}
		&	\leq c \, \smash{h^{\smash{\frac{2}{p'}}}} \,. 
			\label{eq:fem:4error_pres}
	\end{align}
	with a constant $c>0$ depends on the same quantities as $c_1>0$ in Theorem \ref{thm:fem:error_vel_final} and, in addition, on $\rho_{(\phi_{\vert\mathbf{Dv}\vert})\d,\Omega}(\nabla q)$.
\end{proposition}
\begin{proof} 
	\emph{ad \eqref{eq:fem:3error_pres}.} By Lemma \ref{lem:fem:strong_stability} and Lemma \ref{lem:fem:sigma_stability}, we have that\vspace*{-0.5mm}
	\begin{align*}
		\|\mathbf{v}_h\|_{1,d+\epsilon,\Omega} + \|\Sigma_h\mathbf{v}_h\|_{\infty,\Omega} < c_1\,,
	\end{align*}
	which allows to prove an analogue of Lemma \ref{lem:fem:estim_conv_pressure} with $\|\cdot\|_{1,s,\Omega}$ replaced by $\|\cdot\|_{1,2,\Omega}$.
	Hence, we may repeat the proof of Proposition \ref{prop:fem:error_pres}  with $\|\cdot\|_{1,s,\Omega}$ replaced by $\|\cdot\|_{1,2,\Omega}$ to conclude that the claimed \textit{a priori} error estimate for the pressure \eqref{eq:fem:3error_pres} applies.
	
	\emph{ad \eqref{eq:fem:4error_pres}.} First, we note that Lemma \ref{lem:fem:strong_stability} and the embedding $W^{1,d+\epsilon}(\Omega) \hookrightarrow  L^{\infty}(\Omega)$ yield  that $\sup_{h\in (0,1]}{\!\|\mathbf{e}_h\|_{\infty,\Omega}}\!<\!\infty$. Thus,  by \eqref{eq:fem:error_12}, \eqref{cor:fem:error_vel.1}, and real \mbox{interpolation}, we \mbox{get}\vspace*{-0.5mm}
	\begin{equation*}
		\|\mathbf{e}_h\|_{p',\Omega}\leq 
		\left.\begin{cases}
			\|\mathbf{e}_h\|_{p^*,\Omega}^{\theta} \|\mathbf{e}_h\|_{\infty,\Omega}^{1-\theta}
		&\text{ if }p^* < p'\,,\\
			c\,\|\mathbf{e}_h\|_{p^*,\Omega} &\text{ if }p^*\ge p'\,,
		\end{cases} \right\}
		\leq c_1 \smash{h^{\smash{\frac{2}{p'}}}}\,,
	\end{equation*}
	where $\theta = \frac{p^* }{ p'} \geq \frac{2}{p'}$. This allows to prove an analogue of Lemma \ref{lem:fem:estim_conv_pressure} with $\|\cdot\|_{1,s,\Omega}$ replaced by  $\|\cdot\|_{1,p,\Omega}$, which, in turn, allows to prove \eqref{eq:fem:4error_pres} along the lines of Proposition \ref{prop:fem:error_pres}.
\end{proof}

\begin{remark} \label{rem:fem:bootstrap}
	If we had the stronger stability (\textit{cf}.\ Lemma \ref{lem:fem:strong_stability}) at our disposal already in the proof of the error estimates (\textit{cf}.\ Theorem \ref{thm:fem:error_vel_final}), then we could omit the extra regularity assumption on the boundary data $\mathbf{g}_b$ in Assumption \ref{assum:fem:regularity} (\textit{cf}.\ \eqref{eq:Ih1}, \eqref{eq:Ih24}).
	This is why the authors believe that the extra regularity $\mathbf{g}_b \in (W^{2,\ell}(\Omega))^d$ is not necessary. However, it remains an open question how to avoid this assumption in the proof of the error estimates at the first place.
\end{remark}\newpage

\section{Numerical experiments} \label{sec:fem:num_experiments}

\hspace{5mm}In this section, we review the theoretical findings of Section \ref{sec:error_estimates} via numerical experiments. 
Since~in~\cite{JK23_fem_arxiv}, the case $p\hspace{-0.1em}\ge\hspace{-0.1em} \frac{2d}{d+1}$ (\textit{i.e.}, $p\hspace{-0.1em}\ge\hspace{-0.1em} \frac{4}{3}$ if $d\hspace{-0.1em}=\hspace{-0.1em}2$) is already studied, the numerical~\mbox{experiments}~\mbox{focus}~on~the~case~${p\hspace{-0.1em}<\hspace{-0.1em}2}$.\vspace{-2.5mm}

\subsection{Implementation details}\enlargethispage{5mm}

\hspace*{5mm}All experiments were conducted employing the finite element software \texttt{firedrake} (version 0.13.0, \textit{cf}.\ \cite{firedrake}).
In the numerical experiments, 
we deploy the finite element spaces in Remark~\ref{rem:fem:reconstruction} with a discontinuous pressure space. In the case $p<\frac{2d}{d+1}$, the discrete convective term is given via \eqref{eq:fem:defi_conv_recon} and involves a (divergence) reconstruction operator $\Sigma_h\colon X\to X_h+Z_h$ in the sense of Definition~\ref{defi:fem:reconstruction}.~We~incor\-porated the latter not via constructing it explicitly,~but~weakly~instead, \textit{i.e.}, via  extending Problem (\hyperlink{Qh}{Q$_h$}) 
with an unknown $\mathbf{z}_h\in Z_h$ (serving as a placeholder for $\Sigma_h \mathbf{v}_h\in Z_h$) 
in the following ways:

\emph{$\bullet$ P2P0 element/First-order Bernardi--Raugel element.} According to Remark \ref{rem:fem:reconstruction}(i), a divergence reconstruction operator is given by the Fortin interpolation operator ${\Sigma_h\coloneqq \Pi_h^{rt,0}\colon}$ $ X\to Z_h$ (\textit{cf}.\ \cite[Sec.\ 17.1]{EG21}) of the lowest order Raviart--Thomas element $Z_h\coloneqq \mathcal{R}T^0(\mathcal{T}_h)$. We incorporated this interpolation operator into Problem (\hyperlink{Qh}{Q$_h$}) 
via extending this non-linear saddle point problem on the basis of the definition \eqref{eq:fortin_rt.1}. More precisely, we solve the following equivalent augmented problem:

\textsc{Augmented Problem (\hypertarget{AQh0}{AQ$_h^0$}).} Find $(\mathbf{v}_h,\mathbf{z}_h,q_h)\in X_h\times Z_h\times Q_h$ such that $	\mathbf{v}_h = \mathbf{g}_b^h$ in $ \tr X_h$ and 
for every $(\mathbf{w}_h,\mathbf{y}_h,y_h)\in X_h\times Z_h\times Q_h$, there holds
\begin{align*}
	(\mathbf{S}(\mathbf{Dv}_h), \mathbf{Dw}_h)_{\Omega} - (\mathbf{v}_h\otimes  \mathbf{z}_h, \nabla \mathbf{w}_h)_{\Omega} - (q_h, \operatorname{div} \mathbf{w}_h )_{\Omega} &= (\mathbf{f}, \mathbf{w}_h )_{\Omega} \,, \\
	\langle \mathbf{z}_h\cdot \mathbf{n}, \mathbf{y}_h\cdot \mathbf{n}\rangle_{\mathcal{F}_h}&=\langle \mathbf{v}_h\cdot \mathbf{n}, \mathbf{y}_h\cdot \mathbf{n}\rangle_{\mathcal{F}_h}\,,\\
	(\operatorname{div} \mathbf{v}_h, y_h)_{\Omega} &= (g_1^h, y_h)_{\Omega}\,,
\end{align*}
where we exploit that 
$(\mathbf{y}_h\mapsto \mathbf{y}_h\cdot \mathbf{n}|_F)\colon Z_h\to \mathbb{P}^0(F)$ for all $F\in \mathcal{F}_h$ is surjective and we used the notation $\langle \mathbf{z}_h\cdot \mathbf{n}, \mathbf{y}_h\cdot \mathbf{n}\rangle_{\mathcal{F}_h}\coloneqq\sum_{F\in \mathcal{F}_h}{\int_F{(\mathbf{z}_h\cdot \mathbf{n})\,(\mathbf{y}_h\cdot \mathbf{n})\,\mathrm{d}s}}$.

\emph{$\bullet$ Conforming Crouzeix--Raviart element/Second-order Bernardi--Raugel element.} According to Remark \ref{rem:fem:reconstruction}(ii), a divergence reconstruction operator is given via the Fortin interpolation operator $\Sigma_h\coloneqq \Pi_h^{rt,1} \colon X\to Z_h$ (\textit{cf}.\ \cite[Sec.\ 17.1]{EG21}) of the first degree Raviart--Thomas element  $Z_h\coloneqq \mathcal{R}T^1(\mathcal{T}_h)$. We incorporated this operator into  Problem (\hyperlink{Qh}{Q$_h$})  via  extending  this  non-linear  saddle  point  problem based on the definitions \eqref{eq:fortin_rt.1}, \eqref{eq:fortin_rt.2}.
More precisely, we solve the following~equivalent~\mbox{augmented}~\mbox{problem:}

\textsc{Augmented Problem (\hypertarget{AQh1}{AQ$_h^1$}).} Find $(\mathbf{v}_h,\mathbf{z}_h,q_h) \in  X_h\times Z_h\times Q_h$ such that ${\mathbf{v}_h  = \mathbf{g}_b^h}$ in $\tr X_h $ and
for every $(\mathbf{w}_h,\mathbf{y}_h^1,\mathbf{y}_h^0,y_h) \in  V_h\times Z_h|_{\cup\mathcal{F}_h}\times (\mathbb{P}^0(\mathcal{T}_h))^d \times Y_h$, there holds
\begin{align*}
	(\mathbf{S}(\mathbf{Dv}_h),\mathbf{Dw}_h)_{\Omega} - (\mathbf{v}_h\otimes  \mathbf{z}_h, \nabla \mathbf{w}_h)_{\Omega} - (q_h, \operatorname{div} \mathbf{w}_h )_{\Omega} &= (\mathbf{f}, \mathbf{w}_h )_{\Omega} \,, \\
	\langle \mathbf{z}_h\cdot \mathbf{n}, \mathbf{y}_h^1\cdot \mathbf{n}\rangle_{\mathcal{F}_h}&=\langle \mathbf{v}_h\cdot \mathbf{n}, \mathbf{y}_h^1\cdot \mathbf{n}\rangle_{\mathcal{F}_h}\,,\\
	(\mathbf{z}_h, \mathbf{y}_h^0)_{\Omega}&=(\mathbf{v}_h, \mathbf{y}_h^0)_{\Omega}\,,\\
	(\operatorname{div} \mathbf{v}_h, y_h)_{\Omega} &= (g_1^h, y_h)_{\Omega}\,,
\end{align*}
where we exploit that 
$(\mathbf{y}_h\mapsto \mathbf{y}_h\cdot \mathbf{n}|_F)\colon Z_h|_F\to \mathbb{P}^1(F)$ for all $F\in \mathcal{F}_h$ is surjective.

We emphasize the Augmented Problems (\hyperlink{AQh0}{AQ$_h^0$}) and (\hyperlink{AQh1}{AQ$_h^1$}) are designed in such a way that both the formulas \eqref{eq:fortin_rt.1}, \eqref{eq:fortin_rt.2} (where \eqref{eq:fortin_rt.2} is only relevant for Augmented Problem (\hyperlink{AQh1}{AQ$_h^1$})) for the reconstruction operator and 
the structural connection between the discrete pressure space $Y_h$ and the
discrete velocity space $X_h$ are guaranteed, where $\Sigma_h\colon \hspace{-0.1em}X_h \hspace{-0.1em}\to\hspace{-0.1em} X_h + Z_h$ with $\operatorname{div}(Z_h) = Y_h$~(\textit{cf}.~Remark~\ref{rem:fem:reconstruction}).

As discretization of the boundary data, we employ the nodal interpolation of the~trace~lift~${\mathbf{g}_b\in X^s}$, \textit{i.e.}, we set $\mathbf{g}_b^h\coloneqq \Pi_{h,N}^X\mathbf{g}_b\in X_h$, 
 which, according to Lemma \ref{lem:fem:proj_x2_nodal}, given sufficient regularity~of~the~trace~lift $\mathbf{g}_b\in X^s$, is justified. In addition,
since we will restrict to the~case~${g_1=0}$,~we~set~${g_1^h\coloneqq \langle \operatorname{div} \mathbf{g}_b^h\rangle_{\Omega}\in Y^{s'}}$.

We approximate the discrete solutions $(\mathbf{v}_h,\mathbf{z}_h,q_h)\in X_h\times Z_h\times Q_h$ of the Augmented Problems (\hyperlink{AQh}{AQ$_h^0$}) and (\hyperlink{AQh}{AQ$_h^1$}), respectively, using a standard Newton iteration which is deemed to have converged when the Euclidean norm of the residual falls below the tolerance $\tau_{abs} \coloneqq 1\textrm{e}{-}8$.
The linear system emerging in each Newton iteration is solved using a sparse direct solver from \texttt{MUMPS} (version~5.5.0,~\textit{cf}.~\cite{mumps}).\newpage 

\subsection{Experimental setup}\enlargethispage{3mm}\vspace{-0.5mm}

\hspace*{5mm}We use the Augmented Problems (\hyperlink{AQh0}{AQ$_h^0$}) and  (\hyperlink{AQh1}{AQ$_h^1$}), respectively,  (or equivalently the Problems (\hyperlink{Qh}{Q$_h$}),(\hyperlink{Ph}{P$_h$})) to approximate the non-linear system \eqref{eq:fem:main_problem1}--\eqref{eq:fem:main_problem3} on $\Omega=(0,1)^2$ with  $\mathbf{S}\colon   \mathbb{R}^{2\times 2}\to \mathbb{R}^{2\times 2}_{\mathrm{sym}}$, for  every $\mathbf{A}\in\mathbb{R}^{2\times 2}$ defined by\vspace{-0.5mm}
\begin{align*}
	\mathbf{S}(\mathbf{A}) \coloneqq \nu_0\,(\delta+\vert \mathbf{A}^{\textup{sym}}\vert)^{p-2}\mathbf{A}^{\textup{sym}}\,,
\end{align*}  
where  $\nu_0 =100$, $\delta\coloneqq1\textrm{e}{-}5$, and $ p\in (1, 1.5]$.

As manufactured solutions serve the vector field $\mathbf{v}\in X$ and the function $q \in Q^s$, for every $x=(x_1,x_2)\in \Omega $, defined by\vspace{-0.5mm}
\begin{align*}
	\mathbf{v}(x)\coloneqq 	\vert x\vert^\beta(-x_2,x_1)\,,
	\qquad q(x)\vcentcolon = \vert x\vert^{\gamma}-\langle\,\vert \cdot\vert^{\gamma}\,\rangle_\Omega\,,
\end{align*}
\textit{i.e.}, we choose the right-hand side $\mathbf{f}\in  \smash{(L^{p'}(\Omega))^2}$, the divergence  $g_1=0\in L^s(\Omega)$, and boundary data $\mathbf{g}_2\in  \textcolor{black}{(W^{\smash{1-\frac{1}{s},s}}(\partial\Omega))^2}$ accordingly.

For the regularity of the velocity vector field and kinematic pressure, we choose $\beta =  0.01$ and $\gamma= 1-\smash{\frac{2}{p'}}+\beta$, which just yields that Assumption \ref{assum:fem:regularity} is satisfied by $\mathbf{v}$ and $q$. The boundary data $\mathbf{v}|_{ \partial\Omega}$ may, in general, not be regular enough to admit a lift $\mathbf{g}_b \in (W^{2,\ell}(\Omega))^2$,
nevertheless, the EOC presented in Tables \ref{tab:1}--\ref{tab:3} reaches the theoretically
proven rates. This gives another indication that the extra regularity of the boundary data (\textit{cf}.\ Assumption \ref{assum:fem:regularity})~is~probably~not~necessary.

An initial triangulation $\mathcal{T}_{h_0}$, where $h_0=1$, is constructed via subdividing $\overline{\Omega}$
into four triangles along its diagonals. Then,
refined triangulations $\mathcal{T}_{h_i}$, $i=1,\ldots,7$, where ${h_{i+1}=\frac{h_i}{2}}$~for~all~${i=0,\ldots,6}$, are 
obtained by applying uniform refinement; more precisely, the red-refinement rule (\textit{cf}. \cite[Def.~4.8]{Bartels15}).

As  estimation  of  the  convergence  rates,   the  experimental  order  of  convergence (EOC)\vspace{-0.5mm}
\begin{align}
	\texttt{EOC}_i(e_i)\coloneqq\frac{\log(e_{i+1})-\log(e_i)}{\log(h_{i+1})-\log(h_i)}\,, \quad i=1,\dots,6\,,\label{eoc}
\end{align}
is used, where, for every $i=1,\ldots,7$, we denote by $e_i$ a general error quantity.\vspace{-1mm}

\subsection{Quasi-optimality of the error decay rates derived in Theorem \ref{thm:fem:error_vel_final}}\vspace{-0.5mm}

\hspace*{5mm}For $p \in  \{1.1,1.2,1.3, 4/3, 1.4, 1.5\}$, the first-order Bernardi--Raugel element (\textit{cf}.\ Remark~\ref{rem:fem:reconstruction}(i.a)) and the conforming Crouzeix--Raviart element (\textit{cf}.\ Remark \ref{rem:fem:reconstruction}(ii.a)),  for $i=1,\ldots,7$, we compute the error quantities\vspace{-1mm}
\begin{align*} 
		\smash{e_{\mathbf{v},i}^{\mathbf{F}}}\coloneqq\|\mathbf{F}(\mathbf{D}\mathbf{v}_{h_i})-\mathbf{F}(\mathbf{D}\mathbf{v})\|_{2,\Omega}\,,
\end{align*}
and 
the corresponding EOCs, which are presented in Table \ref{tab:1}:

For both the first-order Bernardi--Raugel element and the conforming Crouzeix--Raviart element,     we report
the expected convergence rate of $ \texttt{EOC}_i(\smash{e_{\mathbf{v},i}^{\mathbf{F}}}) \approx 1$, $i=1,\dots,6$, which
confirms the quasi-optimality of the \textit{a priori} error estimates for the velocity vector field derived in Theorem \ref{thm:fem:error_vel_final}. 
For the P2P0 element (\textit{cf}.\ Remark \ref{rem:fem:reconstruction}(i.b)), we observed the same error decay rates as for the first-order Bernardi--Raugel element.\vspace{-0.5mm}

\begin{table}[H]
	\setlength\tabcolsep{4.5pt}
	\centering
	\begin{tabular}{c |c|c|c|c|c|c|c|c|c|c|c|c|c|c|} 
		\cmidrule[\heavyrulewidth](){2-13}
		
		& \multicolumn{6}{c||}{\cellcolor{lightgray}first-order Bernardi--Raugel}   & \multicolumn{6}{c|}{\cellcolor{lightgray}conforming Crouzeix--Raviart}\\ 
		\cmidrule(){2-13}
		
		& \multicolumn{3}{c|}{\cellcolor{lightgray}using \eqref{eq:fem:defi_conv_recon}} & \multicolumn{3}{c||}{\cellcolor{lightgray}using \eqref{eq:fem:defi_conv_temam}}
		& \multicolumn{3}{c|}{\cellcolor{lightgray}using \eqref{eq:fem:defi_conv_recon}} & \multicolumn{3}{c|}{\cellcolor{lightgray}using \eqref{eq:fem:defi_conv_temam}}\\
		\hline 
		
		\multicolumn{1}{|c||}{\cellcolor{lightgray}\diagbox[height=1.1\line,width=0.11\dimexpr\linewidth]{\vspace{-0.6mm}$i$}{\\[-5mm] $p$}}
		& \cellcolor{lightgray}1.1 & \cellcolor{lightgray}1.2 & \cellcolor{lightgray}1.3 & \cellcolor{lightgray}4/3  & \cellcolor{lightgray}1.4  &  \multicolumn{1}{c||}{\cellcolor{lightgray}1.5} &
	\multicolumn{1}{c|}{\cellcolor{lightgray}1.1}  & \cellcolor{lightgray}1.2 &\cellcolor{lightgray}1.3 & \cellcolor{lightgray}4/3 &\cellcolor{lightgray}1.4 & \cellcolor{lightgray}1.5   \\ \toprule\toprule
	\multicolumn{1}{|c||}{\cellcolor{lightgray}$1$}            & 1.012 & 1.011 & 1.009 & 1.009 & 1.008 & \multicolumn{1}{c||}{1.008} & \multicolumn{1}{c|}{1.001} & 1.002 & 1.002 & 1.002 & 1.002 & 1.002 \\ \hline
	\multicolumn{1}{|c||}{\cellcolor{lightgray}$2$}            & 1.010 & 1.010 & 1.010 & 1.010 & 1.010 & \multicolumn{1}{c||}{1.010} & \multicolumn{1}{c|}{1.009} & 1.010 & 1.010 & 1.010 & 1.010 & 1.010 \\ \hline
	\multicolumn{1}{|c||}{\cellcolor{lightgray}$3$}            & 1.008 & 1.008 & 1.009 & 1.009 & 1.009 & \multicolumn{1}{c||}{1.009} & \multicolumn{1}{c|}{1.007} & 1.007 & 1.008 & 1.008 & 1.008 & 1.008 \\ \hline
	\multicolumn{1}{|c||}{\cellcolor{lightgray}$4$}            & 1.007 & 1.007 & 1.008 & 1.008 & 1.008 & \multicolumn{1}{c||}{1.008} & \multicolumn{1}{c|}{1.006} & 1.006 & 1.007 & 1.007 & 1.007 & 1.008 \\ \hline
	\multicolumn{1}{|c||}{\cellcolor{lightgray}$5$}            & 1.006 & 1.007 & 1.007 & 1.007 & 1.007 & \multicolumn{1}{c||}{1.008} & \multicolumn{1}{c|}{1.006} & 1.006 & 1.007 & 1.007 & 1.007 & 1.008 \\ \hline
	\multicolumn{1}{|c||}{\cellcolor{lightgray}$6$}            & 1.006 & 1.006 & 1.007 & 1.007 & 1.007 & \multicolumn{1}{c||}{1.008} & \multicolumn{1}{c|}{1.006} & 1.006 & 1.007 & 1.007 & 1.007 & 1.008 \\ \toprule\toprule
	\multicolumn{1}{|c||}{\cellcolor{lightgray}\small theory}  & 1.000 & 1.000 & 1.000 & 1.000 & 1.000 & \multicolumn{1}{c||}{1.000} & \multicolumn{1}{c|}{1.000} & 1.000 & 1.000 & 1.000 & 1.000 & 1.000 \\ \toprule
	\end{tabular}\vspace{-2mm}
	\caption{Experimental order of convergence: $\texttt{EOC}_i(\smash{e_{\mathbf{v},i}^{\mathbf{F}}})$, ${i=1,\dots,6}$.}
	\label{tab:1}
\end{table}

\subsection{Quasi-optimality of the error decay rates proved in Proposition \ref{prop:fem:error_pres2}}

\hspace*{5mm}For $p \in  \{1.1,1.2,1.3, 4/3, 1.4, 1.5\}$, the first-order Bernardi--Raugel element (\textit{cf}.\ Remark~\ref{rem:fem:reconstruction}(i.a)) and the conforming Crouzeix--Raviart element (\textit{cf}.\ Remark \ref{rem:fem:reconstruction}(ii.a)),  we compute the error quantities\vspace*{-0.75mm}
\begin{align*}
	\left.\begin{aligned} 
		\smash{e_{q,i}^{L^{p'}}}&\coloneqq \|q_{h_i}-q\|_{p',\Omega}\,,\\
		\smash{e_{q,i}^{L^2}}&\coloneqq \|q_{h_i}-q\|_{2,\Omega}\,,
	\end{aligned}\quad\right\}\quad i=0,\dots,6\,,
\end{align*}
and 
the corresponding EOCs, which are presented in Table \ref{tab:2} and Table \ref{tab:3}, respectively:

For the first-order Bernardi--Raugel element and the conforming Crouzeix--Raviart element,  we report
the expected convergence rates of $ \texttt{EOC}_i(\smash{e_{q,i}^{L^{p'}}}) \approx \frac{2}{p'}$, $ i=1,\dots,6$, and  $\texttt{EOC}_i(\smash{e_{q,i}^{L^2}}) \approx 1$, ${i=1,\dots,6}$, 
which confirms the quasi-optimality of the \textit{a priori} error estimates for the pressure  in Proposition \ref{prop:fem:error_pres2}. For the P2P0 element (\textit{cf}. Remark \ref{rem:fem:reconstruction}(i.b)), we observed the same error decay rates as for the first-order Bernardi--Raugel element.

\begin{table}[H]
	\setlength\tabcolsep{4.5pt}
	\centering
	\begin{tabular}{c |c|c|c|c|c|c|c|c|c|c|c|c|c|c|} 
		\cmidrule[\heavyrulewidth](){2-13}
		
		& \multicolumn{6}{c||}{\cellcolor{lightgray}first-order Bernardi--Raugel}   & \multicolumn{6}{c|}{\cellcolor{lightgray}conforming Crouzeix--Raviart}\\ 
		\cmidrule(){2-13}
		
		& \multicolumn{3}{c|}{\cellcolor{lightgray}using \eqref{eq:fem:defi_conv_recon}} & \multicolumn{3}{c||}{\cellcolor{lightgray}using \eqref{eq:fem:defi_conv_temam}}
		& \multicolumn{3}{c|}{\cellcolor{lightgray}using \eqref{eq:fem:defi_conv_recon}} & \multicolumn{3}{c|}{\cellcolor{lightgray}using \eqref{eq:fem:defi_conv_temam}}\\
		\hline 
		
		\multicolumn{1}{|c||}{\cellcolor{lightgray}\diagbox[height=1.1\line,width=0.11\dimexpr\linewidth]{\vspace{-0.6mm}$i$}{\\[-5mm] $p$}}
		& \cellcolor{lightgray}1.1 & \cellcolor{lightgray}1.2 & \cellcolor{lightgray}1.3 & \cellcolor{lightgray}4/3  & \cellcolor{lightgray}1.4  &  \multicolumn{1}{c||}{\cellcolor{lightgray}1.5} &
		\multicolumn{1}{c|}{\cellcolor{lightgray}1.1}  & \cellcolor{lightgray}1.2 &\cellcolor{lightgray}1.3 & \cellcolor{lightgray}4/3 &\cellcolor{lightgray}1.4 & \cellcolor{lightgray}1.5   \\ \toprule\toprule
		\multicolumn{1}{|c||}{\cellcolor{lightgray}$1$}            & 0.351 & 0.057 & 0.010 & -0.02 & 0.061 & \multicolumn{1}{c||}{0.260} & \multicolumn{1}{c|}{0.052} & 0.227 & 0.387 & 0.436 & 0.527 & 0.644 \\ \hline
		\multicolumn{1}{|c||}{\cellcolor{lightgray}$2$}            & 0.602 & 0.558 & 0.403 & 0.432 & 0.524 & \multicolumn{1}{c||}{0.668} & \multicolumn{1}{c|}{0.179} & 0.328 & 0.457 & 0.496 & 0.570 & 0.671 \\ \hline
		\multicolumn{1}{|c||}{\cellcolor{lightgray}$3$}            & 0.321 & 0.282 & 0.422 & 0.470 & 0.558 & \multicolumn{1}{c||}{0.676} & \multicolumn{1}{c|}{0.182} & 0.333 & 0.462 & 0.501 & 0.573 & 0.671 \\ \hline
		\multicolumn{1}{|c||}{\cellcolor{lightgray}$4$}            & 0.224 & 0.315 & 0.456 & 0.497 & 0.574 & \multicolumn{1}{c||}{0.678} & \multicolumn{1}{c|}{0.183} & 0.334 & 0.464 & 0.503 & 0.575 & 0.671 \\ \hline
		\multicolumn{1}{|c||}{\cellcolor{lightgray}$5$}            & 0.195 & 0.330 & 0.463 & 0.503 & 0.577 & \multicolumn{1}{c||}{0.677} & \multicolumn{1}{c|}{0.183} & 0.335 & 0.464 & 0.503 & 0.575 & 0.672 \\ \hline
		\multicolumn{1}{|c||}{\cellcolor{lightgray}$6$}            & 0.187 & 0.334 & 0.465 & 0.504 & 0.577 & \multicolumn{1}{c||}{0.675} & \multicolumn{1}{c|}{0.183} & 0.335 & 0.464 & 0.503 & 0.575 & 0.672 \\ \toprule\toprule
		\multicolumn{1}{|c||}{\cellcolor{lightgray}\small theory}  & 0.181 & 0.333 & 0.462 & 0.500 & 0.571 & \multicolumn{1}{c||}{0.667} & \multicolumn{1}{c|}{0.181} & 0.333 & 0.462 & 0.500 & 0.571 & 0.667 \\ \toprule
		\end{tabular}\vspace{-2mm}
	\caption{Experimental order of convergence: $\texttt{EOC}_i(\smash{e_{q,i}^{L^{p'}}})$, ${i=1,\dots,6}$.}
	\label{tab:2}
\end{table}

\begin{table}[H]
	\setlength\tabcolsep{4.5pt}
	\centering
	\begin{tabular}{c |c|c|c|c|c|c|c|c|c|c|c|c|c|c|} 
		\cmidrule[\heavyrulewidth](){2-13}
		
		& \multicolumn{6}{c||}{\cellcolor{lightgray}first-order Bernardi--Raugel}   & \multicolumn{6}{c|}{\cellcolor{lightgray}conforming Crouzeix--Raviart}\\ 
		\cmidrule(){2-13}
		
		& \multicolumn{3}{c|}{\cellcolor{lightgray}using \eqref{eq:fem:defi_conv_recon}} & \multicolumn{3}{c||}{\cellcolor{lightgray}using \eqref{eq:fem:defi_conv_temam}}
		& \multicolumn{3}{c|}{\cellcolor{lightgray}using \eqref{eq:fem:defi_conv_recon}} & \multicolumn{3}{c|}{\cellcolor{lightgray}using \eqref{eq:fem:defi_conv_temam}}\\
		\hline 
		
		\multicolumn{1}{|c||}{\cellcolor{lightgray}\diagbox[height=1.1\line,width=0.11\dimexpr\linewidth]{\vspace{-0.6mm}$i$}{\\[-5mm] $p$}}
		& \cellcolor{lightgray}1.1 & \cellcolor{lightgray}1.2 & \cellcolor{lightgray}1.3 & \cellcolor{lightgray}4/3  & \cellcolor{lightgray}1.4  &  \multicolumn{1}{c||}{\cellcolor{lightgray}1.5} &
		\multicolumn{1}{c|}{\cellcolor{lightgray}1.1}  & \cellcolor{lightgray}1.2 &\cellcolor{lightgray}1.3 & \cellcolor{lightgray}4/3 &\cellcolor{lightgray}1.4 & \cellcolor{lightgray}1.5   \\ \toprule\toprule
		\multicolumn{1}{|c||}{\cellcolor{lightgray}$1$}            & 0.571 & 0.092 & -0.02 & 0.045 & 0.191 & \multicolumn{1}{c||}{0.364} & \multicolumn{1}{c|}{0.899} & 0.905 & 0.918 & 0.923 & 0.933 & 0.946 \\ \hline
		\multicolumn{1}{|c||}{\cellcolor{lightgray}$2$}            & 1.197 & 1.092 & 1.016 & 1.005 & 0.998 & \multicolumn{1}{c||}{1.007} & \multicolumn{1}{c|}{1.016} & 1.011 & 1.009 & 1.009 & 1.010 & 1.011 \\ \hline
		\multicolumn{1}{|c||}{\cellcolor{lightgray}$3$}            & 1.129 & 1.047 & 1.013 & 1.009 & 1.009 & \multicolumn{1}{c||}{1.014} & \multicolumn{1}{c|}{1.002} & 1.000 & 1.001 & 1.001 & 1.002 & 1.004 \\ \hline
		\multicolumn{1}{|c||}{\cellcolor{lightgray}$4$}            & 1.059 & 1.019 & 1.007 & 1.006 & 1.006 & \multicolumn{1}{c||}{1.008} & \multicolumn{1}{c|}{1.001} & 1.000 & 1.001 & 1.001 & 1.002 & 1.003 \\ \hline
		\multicolumn{1}{|c||}{\cellcolor{lightgray}$5$}            & 1.029 & 1.008 & 1.004 & 1.003 & 1.004 & \multicolumn{1}{c||}{1.005} & \multicolumn{1}{c|}{1.001} & 1.001 & 1.002 & 1.002 & 1.003 & 1.004 \\ \hline
		\multicolumn{1}{|c||}{\cellcolor{lightgray}$6$}            & 1.015 & 1.004 & 1.003 & 1.003 & 1.003 & \multicolumn{1}{c||}{1.004} & \multicolumn{1}{c|}{1.001} & 1.002 & 1.002 & 1.003 & 1.003 & 1.004 \\ \toprule\toprule
		\multicolumn{1}{|c||}{\cellcolor{lightgray}\small theory}  & 1.000 & 1.000 & 1.000 & 1.000 & 1.000 & \multicolumn{1}{c||}{1.000} & \multicolumn{1}{c|}{1.000} & 1.000 & 1.000 & 1.000 & 1.000 & 1.000 \\ \toprule
	\end{tabular}\vspace{-2mm}
	\caption{Experimental order of convergence: $\texttt{EOC}_i(\smash{e_{q,i}^{L^2}})$, ${i=1,\dots,6}$.}
	\label{tab:3}
\end{table}

	\begin{remark}
		As the Tables \ref{tab:2} and \ref{tab:3} reveal, the decay of the pressure error is really low in the first refinement step, for which we suspect the following cause:
		Due to the definition of the analytical pressure $q$, the pressure error is concentrated around the origin. In the first refinement step, the triangulation does not change much near the origin. Especially for piecewise constant discrete pressure (in the first-order Bernardi--Raugel case), this may cause the error decay~to~be~really~low.
\end{remark}

{\setlength{\bibsep}{0pt plus 0.0ex}\small
		
		\bibliographystyle{aomplain}
		\bibliography{references} 
		
}
\end{document}